\pdfoutput=1
\documentclass{article}
\usepackage{spconf}
\usepackage{hyperref,algorithm}
\usepackage{algpseudocode}
\usepackage{mathtools}
\usepackage{balance}
\usepackage[margin=1in]{geometry}             % See geometry.pdf to learn the layout options. There are lots.
\usepackage{graphicx}
\usepackage{amssymb}
\usepackage{amsmath}
\usepackage{amsthm}
\usepackage{amsfonts}
\usepackage{epstopdf}
\usepackage{dblfloatfix}
\usepackage{color}
\usepackage{array}
\usepackage{multirow}
\usepackage{booktabs}
\usepackage{balance}
\usepackage{setspace}
\newcommand{\R}{\mathbb{R}}
\renewcommand{\l}{\ell}
\newcommand{\Identity}{1\!\!1}

\title{SPECTRUM-ADAPTED POLYNOMIAL APPROXIMATION FOR MATRIX FUNCTIONS}

\name{Li Fan$^{\dagger}$ \qquad  David I Shuman$^{\dagger}$ \qquad  Shashanka Ubaru$^{\ddagger}$ \qquad  Yousef Saad$^{\mathsection}$ 
\thanks{Authors' email addresses: \{lfan,dshuman1\}@macalester.edu, Shashanka.Ubaru@ibm.com, saad@umn.edu}
\thanks{MATLAB code for all numerical experiments in this paper is available at \url{http://www.macalester.edu/\textasciitilde dshuman1/publications.html}. It leverages the open access GSPBox \cite{gspbox}.
}}

\address{$^{\dagger}$ Macalester College, Dept. of Mathematics, Statistics, and Computer Science, St. Paul, MN 55105 \\      $^{\ddagger}$ IBM T.J. Watson Research Center, Yorktown Heights, NY 10598  \\
        $^{\mathsection}$ University of Minnesota, Dept. of Computer Science and Engineering, Minneapolis, MN 55455}
    
\begin{document}
\ninept
\maketitle
\begin{abstract}
%Computations of the form $f({\bf A}){\bf b}$ play an important role in applications in signal processing, machine learning, applied mathematics, and other disciplines. 
We propose and investigate two new methods to approximate $f({\bf A}){\bf b}$ for large, sparse, Hermitian matrices ${\bf A}$. The main idea behind both methods is to first estimate the spectral density of ${\bf A}$, and then find polynomials of a fixed order that better approximate the function $f$ on areas of the spectrum with a higher density of eigenvalues. %We compare the proposed methods to  on matrices with a variety of spectral distributions. 
Compared to state-of-the-art methods such as the Lanczos method and truncated Chebyshev expansion, the proposed methods tend to provide more accurate approximations of $f({\bf A}){\bf b}$ at lower polynomial orders, and for matrices ${\bf A}$ with a large number of distinct interior eigenvalues and a small spectral width. 
%Unlike the Lanczos method, the proposed methods yield an approximation of $f$ independent of the choice of $\bf b$, and are amenable to distributed computation. 
\end{abstract}
\section{INTRODUCTION} \label{Se:intro}
%Efficiently computing %Computations of the form 
%$f({\bf A}){\bf b}$, a function of a large, sparse matrix times a vector, is an important component in 
%%play a critical role in 
%applications with high-dimensional data in signal processing, machine learning, biology, physics, chemistry, information theory, and other disciplines.

Efficiently computing $f({\bf A}){\bf b}$, a function of a large, sparse Hermitian matrix times a vector, is an important component in numerous signal processing, machine learning, applied mathematics, and computer science tasks. Application examples include graph-based semi-supervised learning methods \cite{smola}\nocite{belkin_matveeva}-\cite{zhou_bousquet}; graph spectral filtering in graph signal processing \cite{shuman2013emerging}; convolutional neural networks / deep learning \cite{defferrard2016convolutional,bronstein2017geometric}; clustering \cite{tremblay2016compressive,orecchia2012approximating}; approximating the spectral density of a large matrix \cite{lin_spectral_density}; estimating the numerical rank of a matrix \cite{ubaru2016,ubaru2017}; approximating spectral sums such as the log-determinant of a matrix \cite{ubaru2017fast} or the trace of a matrix inverse for applications in physics, biology, information theory, and other disciplines \cite{han}; solving semidefinite programs \cite{arora2007combinatorial}; simulating random walks \cite[Chapter 8]{sachdeva2014faster}; and solving ordinary and partial differential equations \cite{hochbruck}\nocite{friesner}-\cite{gallopoulos}.

%\medskip
%
%\noindent Application examples:
%\begin{itemize}
%\item Filtering in graph signal processing \cite{shuman2013emerging}
%\item Deep learning
%\item Graph-based semi-supervised learning
%\item Approximate spectral densities of large matrices \cite{lin}
%\item Estimate the numerical rank of large matrices \cite{ubaru2016}, \cite{ubaru2017}
%\item Solve ordinary and partial differential equations \cite{hochbruck}\nocite{friesner}-\cite{gallopoulos}
%\item Approximate spectral sums such as the log-determinant of a large matrix \cite{ubaru2017fast} or the trace of a matrix inverse for applications in physics, biology, information theory, and other disciplines \cite{han}
%\end{itemize}

References \cite[Chapter 13]{higham}, \cite{davies2005computing}\nocite{frommer}-\cite{moler2003nineteen} survey different approaches to this well-studied problem of efficiently computing 
\begin{align}\label{Eq:problem}
f({\bf A}){\bf b}:={\bf V}f({\boldsymbol \Lambda}){\bf V}^{\top}{\bf b},
\end{align}
where the columns of ${\bf V}$ are the eigenvectors of the Hermitian matrix ${\bf A}\in \R^{N \times N}$; ${\boldsymbol \Lambda}$ is a diagonal matrix whose diagonal elements are the corresponding eigenvalues of ${\bf A}$, which we denote by $\lambda_1, \lambda_2, \ldots, \lambda_N$; and $f({\boldsymbol \Lambda})$ is a diagonal matrix whose $k$th diagonal entry is given by $f(\lambda_k)$. For large matrices, it is not practical to explicitly compute the eigenvalues of ${\bf A}$ in order to approximate \eqref{Eq:problem}. Rather, the most common %approximation 
techniques, all of which avoid a full eigendecomposition of ${\bf A}$, include (i) truncated orthogonal polynomial expansions, including Chebyshev \cite{druskin}\nocite{saad2006filtered}-\cite{chen_saad} and Jacobi;
%Laguerre \cite{sheehan2010computing}; 
(ii) rational approximations \cite[Section 3.4]{frommer}; (iii) Krylov subspace methods such as the Lanczos method \cite{druskin}, \cite{druskin1998extended}\nocite{eiermann2006restarted,afanasjew2008implementation}-\cite{frommer2017radau}; and (iv) quadrature/contour integral methods \cite[Section 13.3]{higham}. 
%{\color{red}
%%\noindent Existing methods, include links to surveys:
%\begin{itemize}
%%\item Polynomial approximations: Chebyshev, Legendre, Laguerre, polynomial approximation of splines
%%\item Rational approximation
%\item Jacobi
%\item Krylov / Lanczos / Conjugate gradient / Algebraic multigrid
%\item Contour integral
%\end{itemize}
%}

Our focus in this work is on polynomial approximation methods. Let $p_K(\lambda)=c_0+\sum_{k=1}^K c_k \lambda^k$ be a degree $K$ polynomial approximation to the function $f$ on a known interval $[\underline{\lambda},\overline{\lambda}]$ containing all of the eigenvalues of ${\bf A}$. Then the approximation $p_K({\bf A}){\bf b}$ can be computed recursively, either through a three-term recurrence for specific types of polynomials (see Section \ref{Se:spectrum_adapted} for more details), or 
through a nested multiplication iteration \cite[Section 9.2.4]{golub}, letting ${\bf x}^{(0)} = c_K {\bf b}$, and then iterating
\begin{align*}%\label{Eq:horner}
{\bf x}^{(l)}=c_{K-l}{\bf b} + {\bf A}{\bf x}^{(l-1)},~l=1,2,\ldots,K.
\end{align*}
The computational cost of either of these approaches is dominated by multiplying the sparse matrix ${\bf A}$ by $K$ different vectors. The approximation error is bounded by 
\begin{align}
||f({\bf A})-p_K({\bf A})||_2&=\max_{\l=1,2,\ldots,N}|f(\lambda_{\l})-p_K(\lambda_{\l})| \label{Eq:exact_error} \\ 
&\leq \sup_{\lambda \in [\underline{\lambda},\overline{\lambda}]}|f(\lambda)-p_K(\lambda)|. \label{Eq:range_error}
\end{align}
If, for example, $p_K$ is a degree $K$ truncated Chebyshev series approximation of an analytic function $f$, the upper bound in \eqref{Eq:range_error} converges geometrically to 0 as $K$ increases, at a rate of ${\mathcal{O}}\left({\rho^{-K}}\right)$, where $\rho$ is the radius of an open Bernstein ellipse on which $f$ is analytic and bounded (see, e.g., \cite[Theorem 5.16]{handscomb}, \cite[Theorem 8.2]{atap}).
%{\color{red} Convergence, Chebyshev example.}
In addition to the computational efficiency and convergence guarantees, a third advantage of polynomial approximation methods is that they can be implemented in a distributed setting \cite{shuman_distributed_SIPN_2018}. A fourth advantage is that the $i$th element of $p_K({\bf A}){\bf b}$ only depends on the elements of ${\bf b}$ within $K$ hops of $i$ on the graph associated with ${\bf A}$. This localization property is important in many graph-based data analysis applications (e.g., graph spectral filtering \cite{LTS-ARTICLE-2009-053}, deep learning \cite{defferrard2016convolutional}).

\begin{figure}[t]
\centerline{\includegraphics[width=7.5cm]{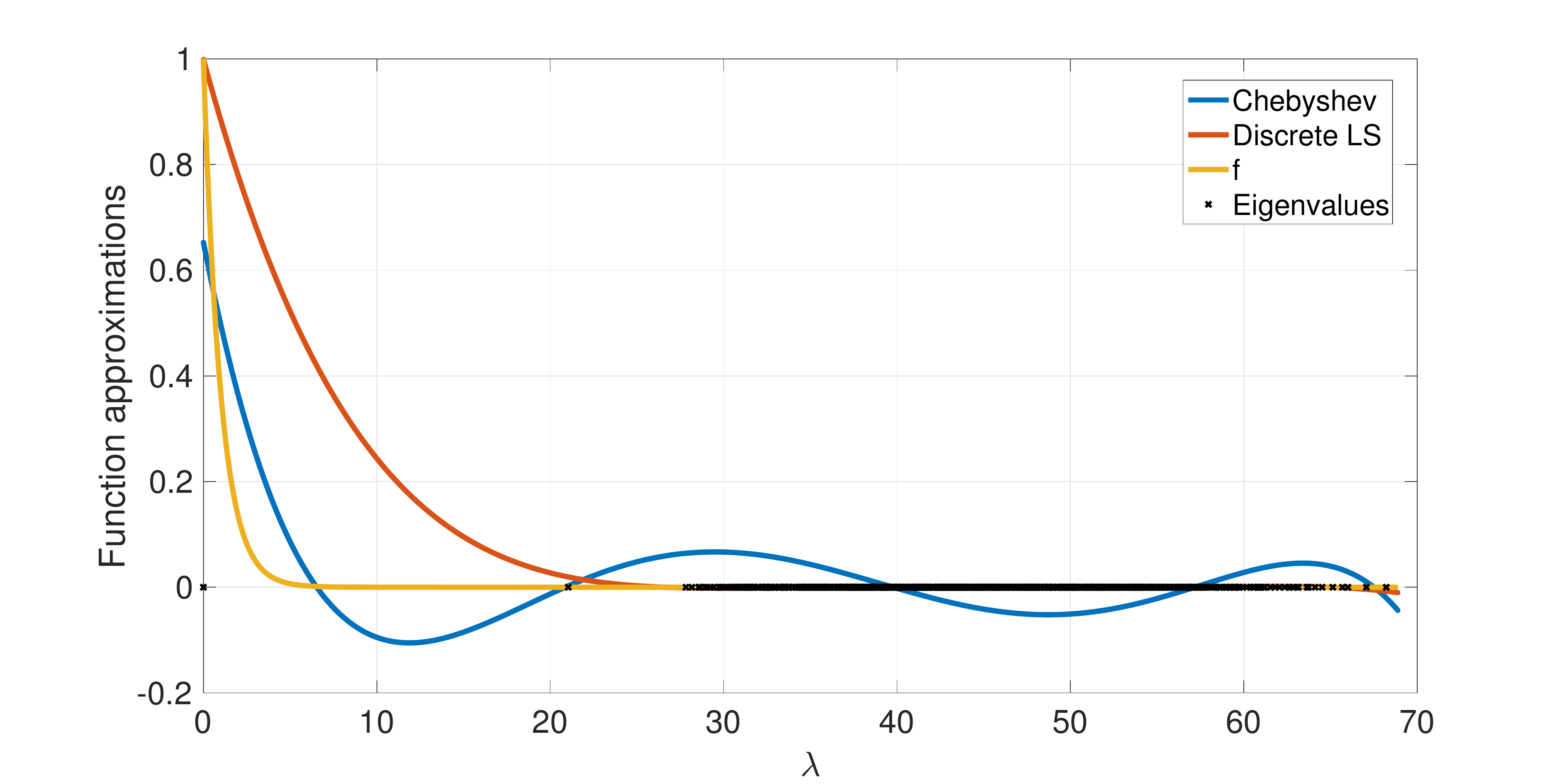}}
\centerline{\includegraphics[width=7.5cm]{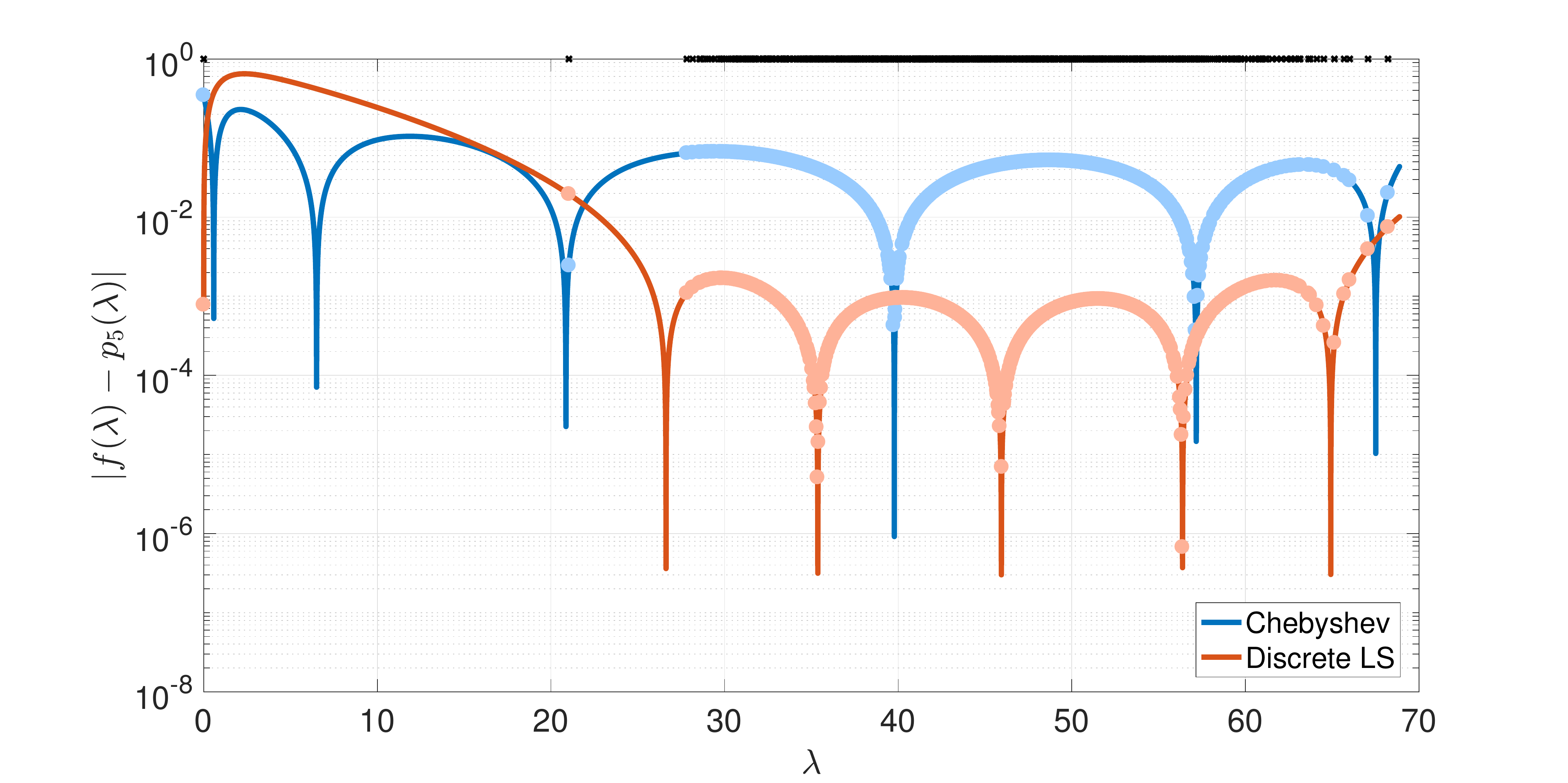}}
\caption{Degree 5 polynomial approximations of the function $f(\lambda)=e^{-\lambda}$ of the graph Laplacian of a random Erd\"{o}s-Renyi graph with 500 vertices and edge probability 0.2. The discrete least squares approximation incurs larger errors in the lower end of the spectrum. However, since the eigenvalues are concentrated at the upper end of the spectrum, it yields a lower approximation error $||f({\bf A})-p_5({\bf A})||_2$.}
\label{Fig:dls}
\end{figure}

While the classical truncated orthogonal polynomial expansion methods (e.g., Chebyshev, Legendre, Jacobi) aim to approximate the function $f$ throughout the full interval $[\underline{\lambda},\overline{\lambda}]$, it is only the polynomial approximation error at the eigenvalues of ${\bf A}$ that affects the overall error in \eqref{Eq:exact_error}. With knowledge of the complete set of eigenvalues, we could do better, for example, by fitting a degree $K$ polynomial via the discrete least squares problem $\min_{p \in {\cal P}_K} \sum_{\l=1}^N \left[f(\lambda_{\l})-p(\lambda_{\l})\right]^2$. In Fig. \ref{Fig:dls}, we show an example of such a discrete least squares fitting. The resulting approximation error $||f({\bf A})-p_K({\bf A})||_2$ for $K=5$ is 0.020, as opposed to 0.347 for the degree 5 truncated Chebyshev approximation. This is despite the fact that $\sup_{\lambda \in [\underline{\lambda},\overline{\lambda}]}|f(\lambda)-p_K(\lambda)|$ is equal to 0.650 for the discrete least squares approximation, as opposed to 0.347 for the Chebyshev approximation.

While in our setting we do not have access to the complete set of eigenvalues, our approach in this work is to leverage recent developments in efficiently estimating the spectral density of the matrix ${\bf A}$, to adapt the polynomial to the spectrum in order to achieve better approximation accuracy at the (unknown) eigenvalues. After reviewing spectral density estimation in the next section, we present two new classes of spectrum-adapted approximation techniques in Section \ref{Se:spectrum_adapted}. We conclude with numerical experiments, and a discussion of the situations in which the proposed methods work better than the state-of-the-art methods.

%%\noindent Motivation and approach:
%\begin{itemize}
%%\item Localization is important in many graph-based applications (filtering, deep learning)
%\item Existing methods try to approximate the function throughout the range of the spectrum; however, all that matters is the approximation at the eigenvalues
%\item Motivating example: Chebyshev vs. discrete least squares with known eigenvalues; include figures of the approximation and errors
%\item For large matrices, we do not want to compute all of these eigenvalues
%\item Our approach here is to leverage recent developments in efficiently estimating the spectral density of the matrix $A$ to adapt the polynomial to the spectrum in order to achieve better approximation accuracy at the (unknown) eigenvalues
%\item Mention that Lanczos is also spectrum-adapted
%\end{itemize}

\section{SPECTRAL DENSITY ESTIMATION}
\begin{figure}[bt] 
\begin{minipage}[m]{0.32\linewidth}
\centerline{\small{~~~~gnp}}
\centerline{~~\includegraphics[width=1.1\linewidth]{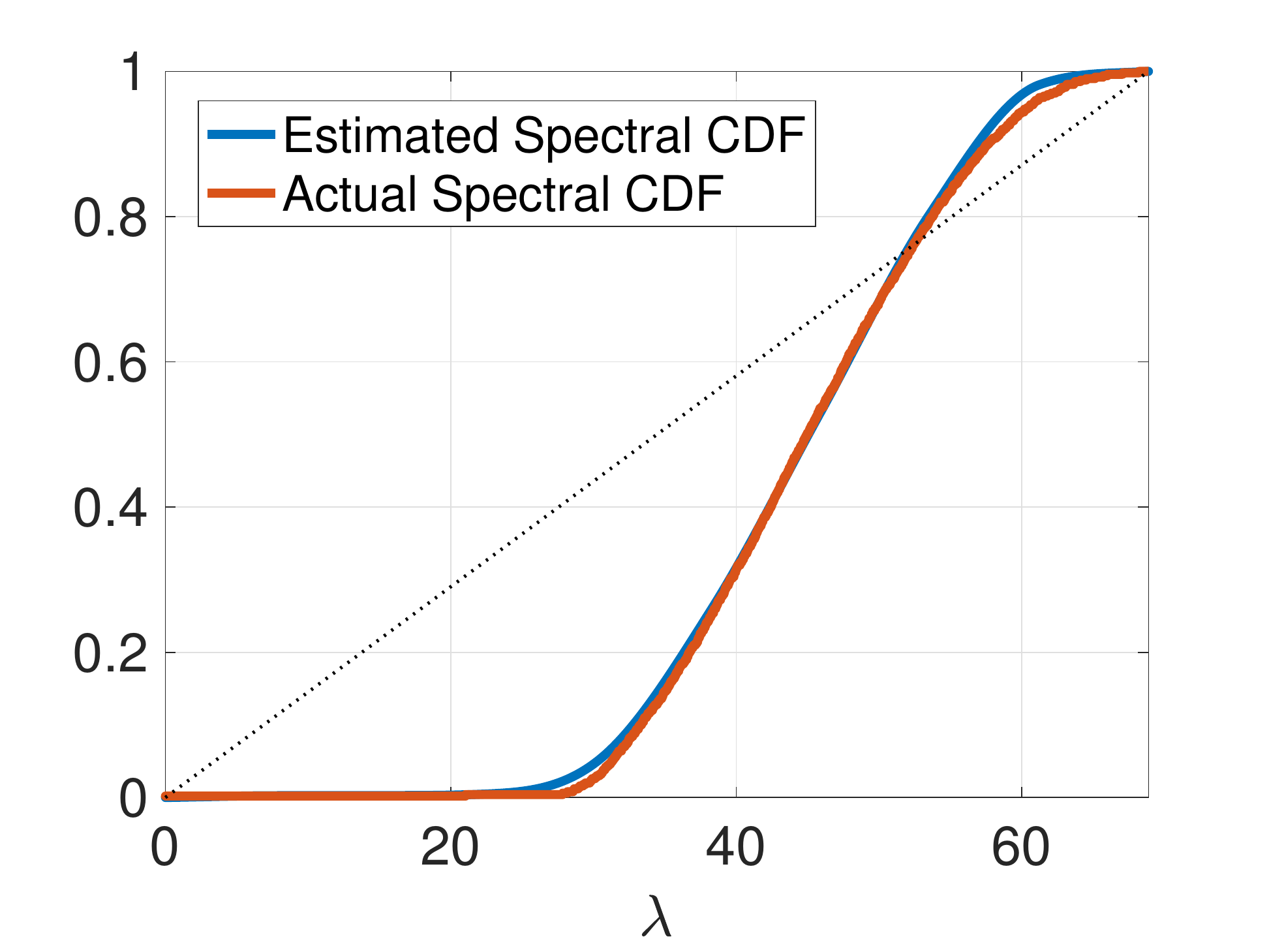}}
\end{minipage}
\begin{minipage}[m]{0.32\linewidth}
\centerline{\small{~~~~minnesota}}
\centerline{~~\includegraphics[width=1.1\linewidth]{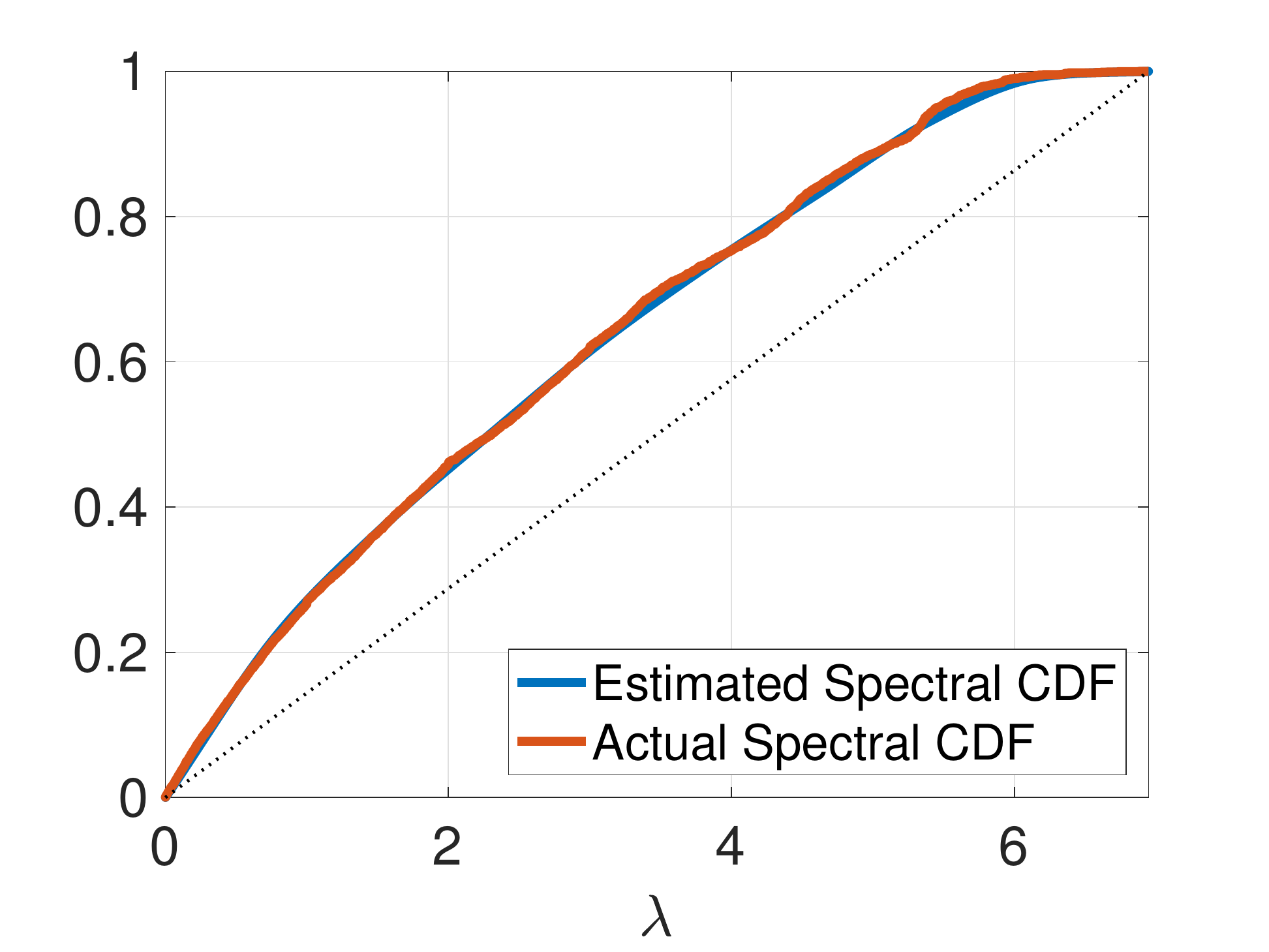}}
\end{minipage}
\begin{minipage}[m]{0.32\linewidth}
\centerline{\small{~~~~net25}}
\centerline{~~\includegraphics[width=1.1\linewidth]{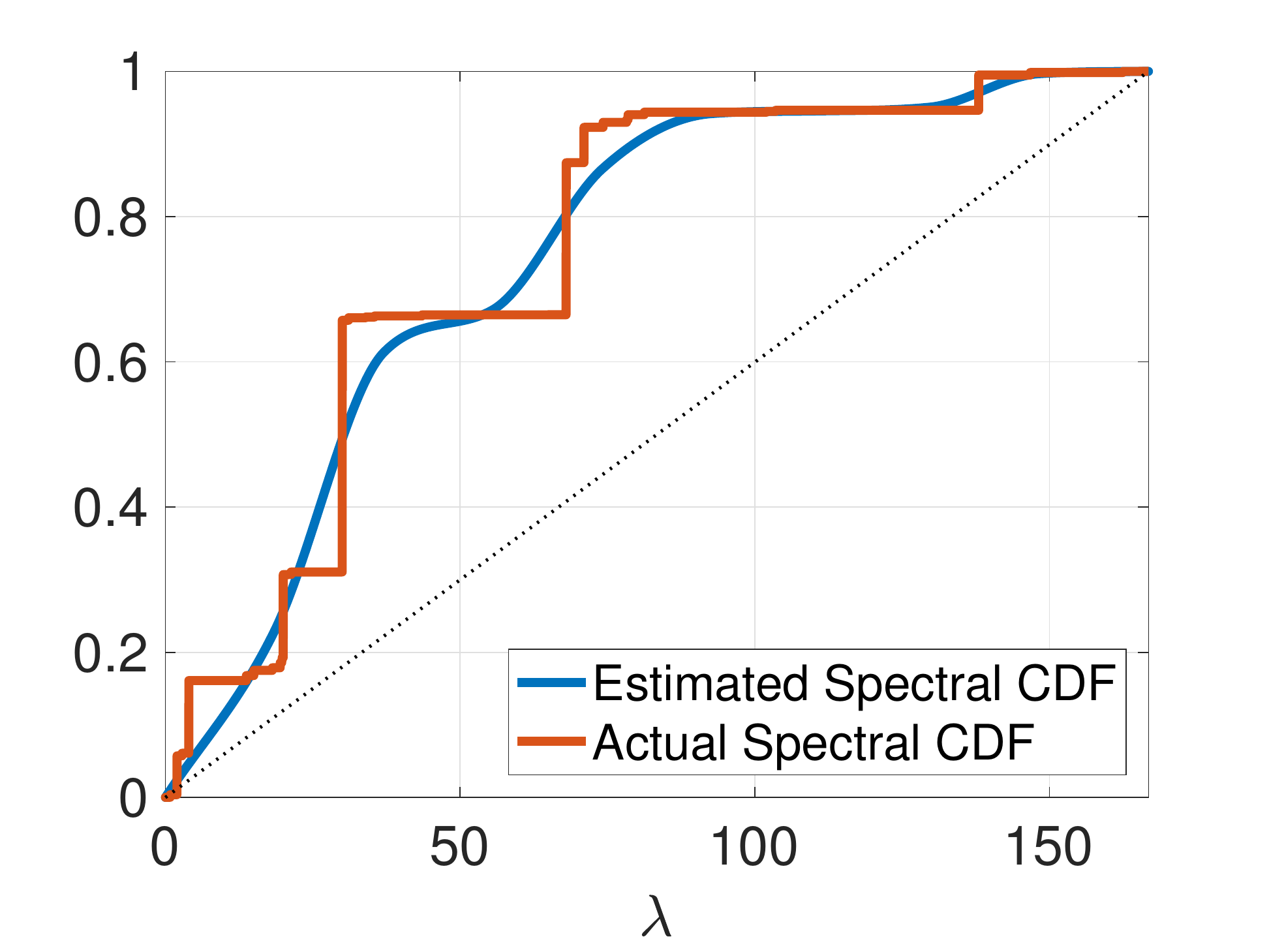}}
\end{minipage}
\\
\vspace{0.03\linewidth}

\begin{minipage}[m]{0.32\linewidth}
\centerline{\small{~~~~si2}}
\centerline{~~\includegraphics[width=1.1\linewidth]{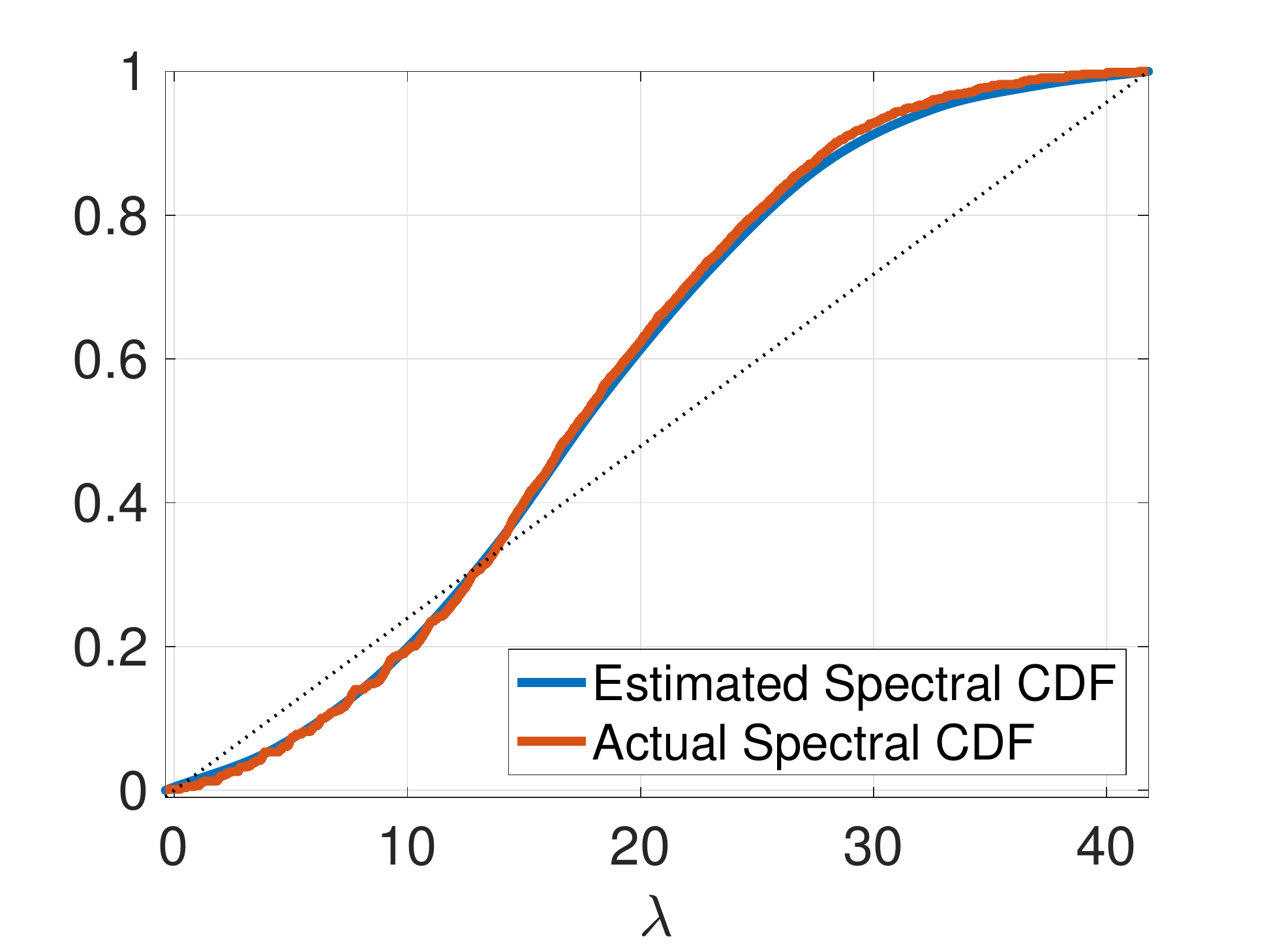}}
\end{minipage} 
\begin{minipage}[m]{0.32\linewidth}
\centerline{\small{~~~~cage9}} 
\centerline{~~\includegraphics[width=1.1\linewidth]{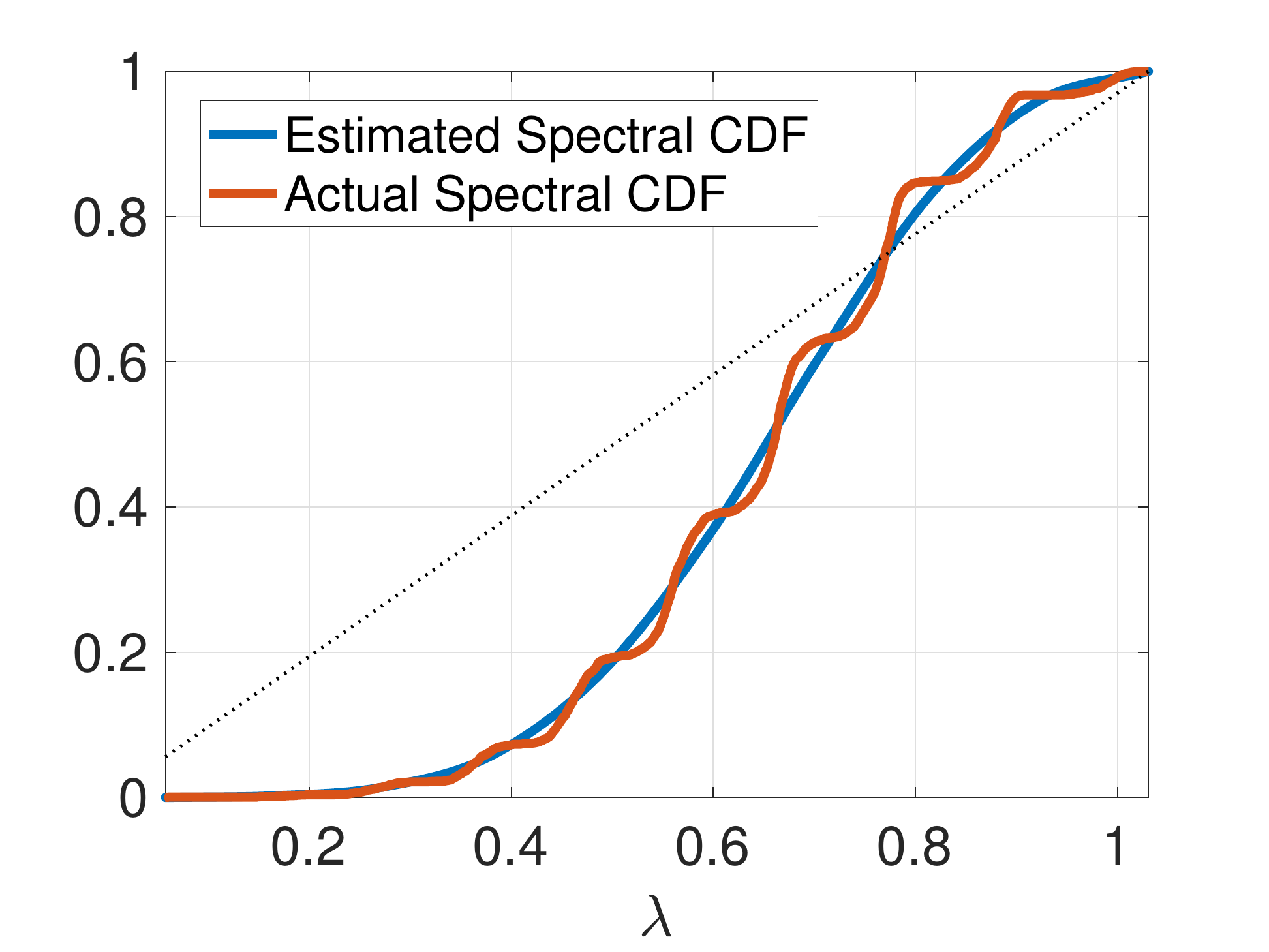}}
\end{minipage}
\begin{minipage}[m]{0.32\linewidth}
\centerline{\small{~~~~saylr4}} 
\centerline{~~\includegraphics[width=1.1\linewidth]{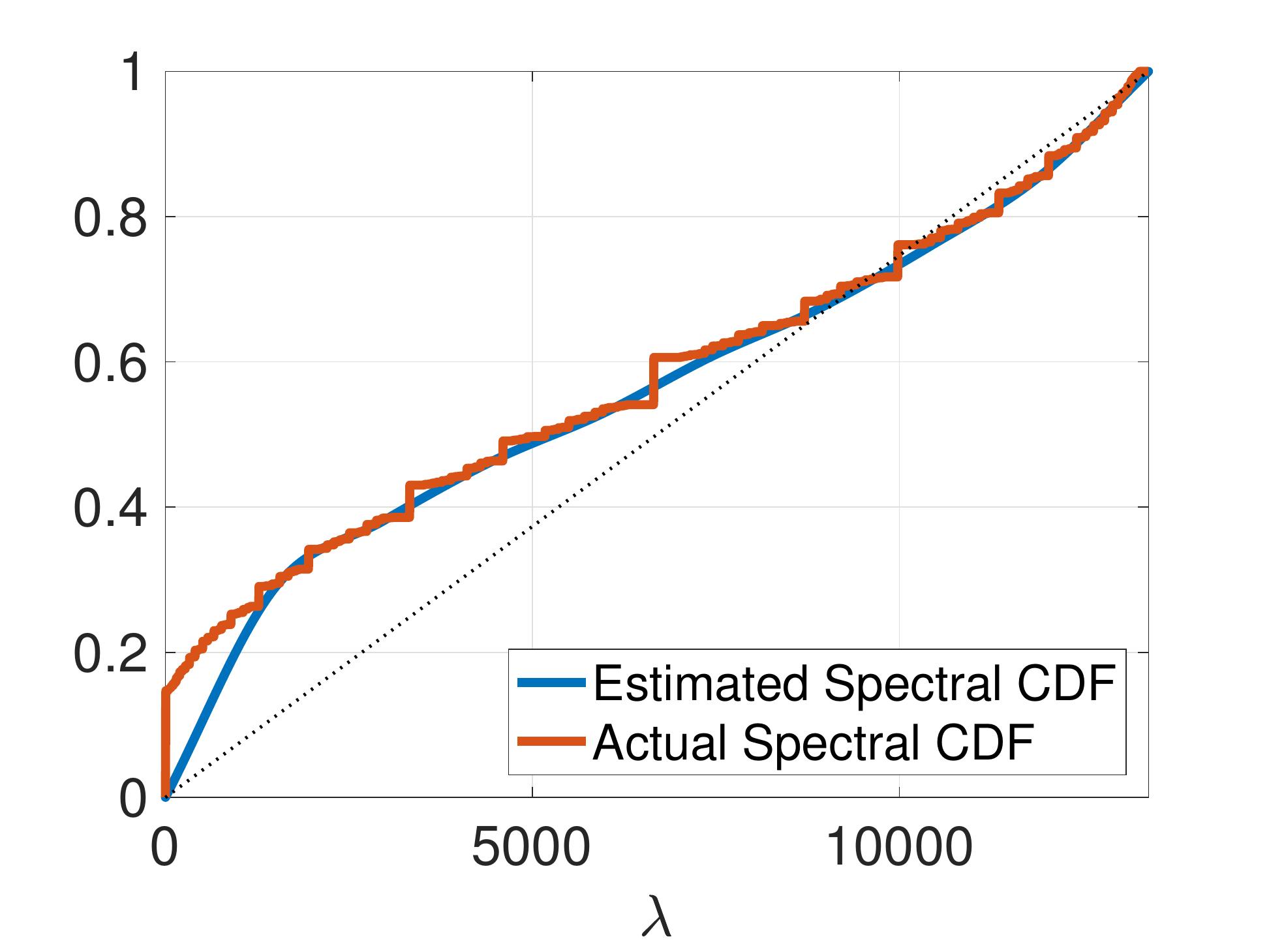}}
\end{minipage}
\caption{Estimated and actual cumulative spectral density functions for six real, symmetric matrices {\bf A}. We estimate the eigenvalue counts for $T=10$ linearly spaced points on $[\underline{\lambda}, \overline{\lambda}]$ via \eqref{Eq:eig_count}, with degree $K_{\Theta}=30$ polynomials and $J=10$ random vectors ${\bf x}^{(j)}$.}\label{Fig:cdf}
\end{figure}

The \emph{cumulative spectral density function} or \emph{empirical spectral cumulative distribution} of the matrix ${\bf A}$ is defined as 
\begin{align} \label{Eq:cdf}
P_{\lambda}(z):=\frac{1}{N}\sum_{\l=1}^{N} \Identity_{\left\{\lambda_{\l}\leq z\right\}},
\end{align}
and the \emph{spectral density function} \cite[Chapter 6]{van_mieghem}) (also called the \emph{Density of States} or \emph{empirical spectral distribution} \cite[Chapter 2.4]{tao_random_matrix}) of ${\bf A}$ is the probability measure defined as 
%\begin{align*}
$p_{\lambda}(z):=\frac{1}{N}\sum_{\l=1}^{N} \Identity_{\left\{\lambda_{\l}=z\right\}}.$
%\end{align*}
Lin et al. \cite{lin_spectral_density} provide an %excellent 
overview of methods to approximate these functions. In this work, we use a variant of the Kernel Polynomial Method (KPM) \cite{silver1994densities}\nocite{silver1996kernel}-\cite{wang1994calculating} described in \cite{lin_spectral_density,mcsfb} to estimate the cumulative spectral density function $P_{\lambda}(z)$ of ${\bf A}$. Namely, for each of $T$ linearly spaced points $\xi_i$ between $\underline{\lambda}$ and $\overline{\lambda}$, we estimate the number of eigenvalues less than or equal to $\xi_i$ via Hutchinson's stochastic trace estimator \cite{hutchinson}:
\begin{align}\label{Eq:eig_count}
\eta_i =\mbox{tr}\bigl(\Theta_{\xi_i}({\bf A})\bigr)
=\mathbb{E}[{\bf x}^{\top}\Theta_{\xi_i}({\bf A}){\bf x}]   
%&\approx \frac{1}{J} \sum_{j=1}^J {{\bf x}^{(j)}}^{\top}\Theta_{\xi_i}(\L){\bf x}^{(j)} \nonumber  \\
&\approx \frac{1}{J} \sum_{j=1}^J {{\bf x}^{(j)}}^{\top}\tilde{\Theta}_{\xi_i}({\bf A}){\bf x}^{(j)},
\end{align}
where each ${\bf x}^{(j)}$ is random vector with each component having an independent and identical standard normal distribution, and $\tilde{\Theta}_{\xi_i}$ is a Jackson-Chebyshev polynomial approximation to $\Theta_{\xi_i}(\lambda):=\Identity_{\left\{\lambda \leq \xi_i\right\}}$ \cite{di2016efficient,puy_structured_sampling}. As in \cite{shuman2013spectrum}, we then form an approximation $\tilde{P}_{\lambda}(z)$ to $P_{\lambda}(z)$ by performing monotonic piecewise cubic interpolation \cite{fritsch} on the series of points $\left\{\left(\xi_i,\frac{\eta_i}{N}\right)\right\}_{i=1,2,\ldots,T}$. Analytically differentiating $\tilde{P}_{\lambda}(z)$ yields an approximation $\tilde{p}_{\lambda}(z)$ to the spectral density function $p_{\lambda}(z)$. Since $\tilde{P}_{\lambda}(z)$ is a monotonic cubic spline, we can also analytically compute its inverse function $\tilde{P}_{\lambda}^{-1}(y)$. Fig. \ref{Fig:cdf} shows examples of the estimated cumulative spectral density functions for six real, symmetric matrices ${\bf A}$: the graph Laplacians of the Erd\"{o}s-Renyi graph (gnp) from Fig. \ref{Fig:dls} and the Minnesota traffic network \cite{gleich} ($N=2642$), and the net25 ($N=9520$), si2 ($N=769$), cage9 ($N=3534$), and saylr4 ($N=3564$) matrices from the  SuiteSparse Matrix Collection \cite{suitesparse}.\footnote{We use $\frac{{\bf A}+{\bf A^{\top}}}{2}$ for cage9, and for net25 and saylr4, we generate graph Laplacians based on the off-diagonal elements of ${\bf A}$.}  The computational complexity of forming the estimate $\tilde{P}_{\lambda}(z)$ is ${\cal O}(MJK_{\Theta})$, where $M$ is the number of nonzero entries in ${\bf A}$, $J$ is the number of random vectors in \eqref{Eq:eig_count} (in our experiments, $J=10$ suffices), and $K_{\Theta}$ is the degree of the Jackson-Chebyshev polynomial approximations $\tilde{\Theta}_{\xi_i}$ \cite{mcsfb}. While this cost is non-negligible if computing $f({\bf A}){\bf b}$ for a single $f$ and a single ${\bf b}$, it only needs to be computed once for each ${\bf A}$ if repeating this calculation for multiple functions $f$ or multiple vectors ${\bf b}$, as is often the case in the applications mentioned above.

%\begin{itemize}
%\item Mention other methods
%\item KPM method
%\item CDF and PDF analytically
%\item Examples of CDFs
%\end{itemize}
%
%
%
%If we place the $J$ random vectors into the columns of an $N \times J$ matrix ${\bf X}$, the computational cost of estimating the spectral distribution is dominated by computing 
%\begin{align} \label{Eq:hutch4}
%\tilde{\Theta}_{\xi_i}(\L){\bf X}=\sum_{k=0}^K \alpha_k \bar{T}_k(\L){\bf X}
%\end{align}
% for each $\xi_i$. Yet, %the $\alpha_k$'s for each ${\Theta}_{\xi_i}$ can be computed in closed form, and 
%we only need to compute $\{ \bar{T}_k(\L){\bf X}\}_{k=0,1,\ldots,K}$ recursively once, as this sequence can be reused for each $\xi_i$, with different choices of the $\alpha_k$'s. Therefore, the overall computational cost is ${\cal O}(KJ|\E|)$.

\section{SPECTRUM-ADAPTED METHODS} \label{Se:spectrum_adapted}
In this section, we introduce two new classes of degree $K$ polynomial approximations $p_K({\bf A}){\bf b}$ to $f({\bf A}){\bf b}$, both of which leverage the estimated cumulative spectral density function $\tilde{P}_{\lambda}(z)$. 

\subsection{Spectrum-adapted polynomial interpolation} \label{Se:interpolation}
\begin{figure}[tbh]
\centerline{\includegraphics[width=4.1cm]{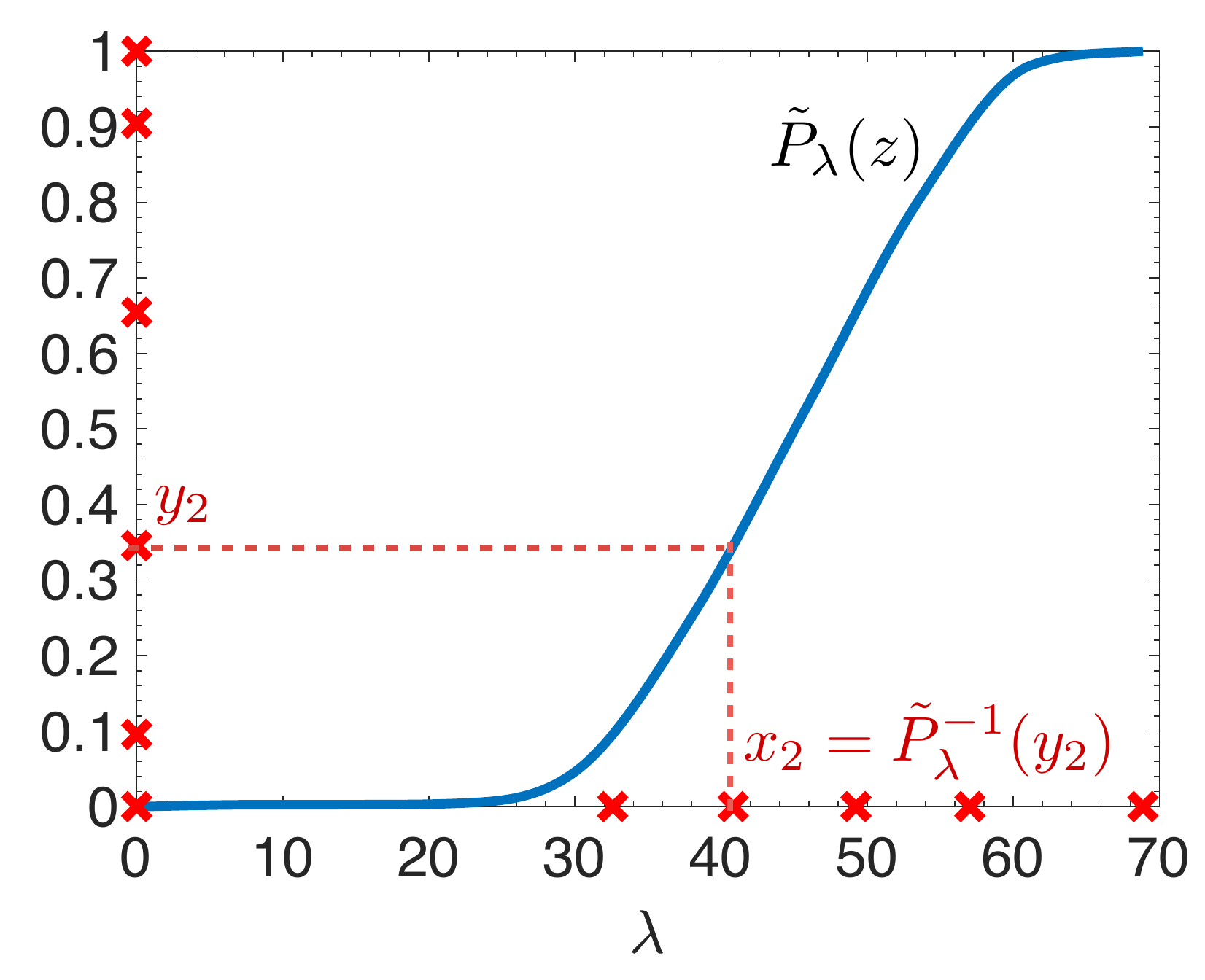}}
\caption{Construction of six interpolation points for the same graph Laplacian matrix described in Fig. \ref{Fig:dls}. The interpolation points $\{x_k\}$ on the horizontal axis are computed by applying the inverse of the estimated cumulative spectral density function to the initial Chebyshev points $\{y_k\}$ on the vertical axis.}
\label{Fig:inverse}
\end{figure}

In the first method, we take $y_k:=\frac{\cos\left(\frac{k\pi}{K}\right)+1}{2}$, for $k=0, 1, \ldots, K$, which are the $K+1$ extrema of the degree $K$ Chebyshev polynomial shifted to the interval $[0,1]$. We then 
warp these points via the inverse of the estimated cumulative spectral density function by setting $x_k=P_{\lambda}^{-1}(y_k)$, %for $k=0, 1, \ldots, K$, 
before finding the unique degree $K$ polynomial interpolation through the points $\{(x_k, f(x_k))\}_{k=0,1,\ldots,K}$. As shown in Fig. \ref{Fig:inverse}, a higher density of the warped points $\{x_k\}$ fall in higher density regions of the spectrum of ${\bf A}$.
\vspace{-.1in}
%\noindent 1. Warp $K$ points with inverse of CDF, interpolate.

\subsection{Spectrum-adapted polynomial regression / orthogonal polynomial expansion} 
A second approach is to solve the weighted least squares polynomial regression problem
\vspace{-.15in}
\begin{align*}
\min_{p \in {\cal P}_K} \sum_{m=1}^M w_m \left[f(x_m)-p(x_m)\right]^2,
\end{align*}
where the abscissae $\{x_m\}_{m=1,2,\ldots,M}$ and weights $\{w_m\}_{m=1,2,\ldots,M}$ are chosen to capture the estimated spectral density function. 
We investigated several methods to set the points (e.g., linearly spaced points, Chebyshev points on the interval $[\underline{\lambda},\overline{\lambda}]$, Chebyshev points on each subinterval $[\xi_i,\xi_{i+1}]$, and warped points via the inverse of the estimated cumulative spectral density function as in Section \ref{Se:interpolation}) and weights (e.g., the analytically computed estimate $\tilde{p}_{\lambda}$ of the spectral density function, a discrete estimate of the spectral density function
 based on the eigenvalue counts in \eqref{Eq:eig_count}, the original KPM density of states method based on a truncated Chebyshev expansion \cite[Eq. 3.11]{lin_spectral_density}, or equal weights %when using 
for warped points). In the numerical experiments, we use $M$ evenly spaced points on the interval $[\underline{\lambda},\overline{\lambda}]$ (i.e., $x_m=\frac{m-1}{M-1}(\overline{\lambda}-\underline{\lambda})+\underline{\lambda}$), and set the weights to be $w_m=\tilde{p}_{\lambda}(x_m)$.
 
An alternative way to view this weighted least squares method \cite{forsythe1957generation} is as a truncated expansion in polynomials orthogonal with respect to the discrete measure $d\lambda_M$ with finite support at the points $\{x_m\}$, and an associated inner product \cite[Section 1.1]{gautschi2004orthogonal} %{\color{red} [use square root of w here?]}
\begin{align*}
\langle f,g \rangle_{d\lambda_M} = \int_{\R} f(t)g(t)d\lambda_M = \sum_{m=1}^M w_m f(x_m) g(x_m).
\end{align*}
The $M$ discrete monic orthogonal polynomials $\{\pi_{k,M}\}_{k=0,1,M-1}$ satisfy the three-term recurrence relation \cite[Section 1.3]{gautschi2004orthogonal}
\begin{align}\label{Eq:three_term}
&\pi_{k+1,M}(x)=(x-\alpha_{k,M})\pi_{k,M}(x)-\beta_{k,M} \pi_{k-1,M}(x), \nonumber \\
&\hspace{1.6in}k=0,1,\ldots,M-1,
\end{align}
with $\pi_{-1,M}(x)=0$, $\pi_{0,M}(x)=1$, $\beta_{0,M}=\sum_{m=1}^M w_m$, 
\begin{align*}
\alpha_{k,M}=\frac{\langle t\pi_{k,M},\pi_{k,M} \rangle_{d\lambda_M}}{\langle \pi_{k,M},\pi_{k,M} \rangle_{d\lambda_M}},~k=0,1,\ldots,M-1,
\end{align*}
\begin{align*}
\hbox{and } \beta_{k,M}=\frac{\langle \pi_{k,M},\pi_{k,M} \rangle_{d\lambda_M}}{\langle \pi_{k-1,M},\pi_{k-1,M} \rangle_{d\lambda_M}},~k=1,2,\ldots,M-1.
\end{align*}
%{\color{blue}
%\begin{itemize}
%\item Three term recurrence relation and discrete orthogonal polynomials
%\end{itemize}
%}
Given the abscissae $\{x_m\}$ and weights $\{w_m\}$, the three-term recursion coefficients $\{\alpha_{k,M}\}_{k=0,1,\ldots,M-1}$ and $\{\beta_{k,M}\}_{k=1,2,\ldots,M-1}$ can also be computed through a stable Lanczos type algorithm on an $(M+1) \times (M+1)$ matrix \cite[Section 2.2.3]{gautschi2004orthogonal}, \cite{gragg1984numerically}. In matrix-vector notation, the vectors ${\boldsymbol \pi}_{k,M} \in \R^{M}$, which are the discrete orthogonal polynomials evaluated at the $M$ abscissae, can be computed iteratively by the relation
\begin{align*}
&{\boldsymbol \pi}_{k+1,M}=(\hbox{diag}(\{x_m\})-\alpha_{k,M}{\bf I}_M){\boldsymbol \pi}_{k,M}-\beta_{k,M} {\boldsymbol \pi}_{k-1,M}, \\
&\hspace{2in}k=0,1,\ldots,M-1,
\end{align*}
with ${\boldsymbol \pi}_{-1,M}={\bf 0}_M$ and ${\boldsymbol \pi}_{0,M}={\bf 1}_M$.
Finally, the degree $K$ polynomial approximation to $f({\bf A}){\bf b}$ is computed as  
\begin{align*}
p_K({\bf A}){\bf b}=\sum_{k=0}^K \frac{\langle f,\pi_{k,M}\rangle_{d\lambda_M}}{\langle \pi_{k,M},\pi_{k,M} \rangle_{d\lambda_M}} \pi_{k,M}({\bf A}){\bf b},
\end{align*}
with $\pi_{-1,M}({\bf A}){\bf b}={\bf 0}_N$, $\pi_{0,M}({\bf A}){\bf b}={\bf b}$, and 
\begin{align*}
&\pi_{k+1,M}({\bf A}){\bf b}=({\bf A}-\alpha_{k,M}{\bf I}_N)\pi_{k,M}({\bf A}){\bf b}-\beta_{k,M} \pi_{k-1,M}({\bf A}){\bf b}, \\
&\hspace{1.1in}k=0,1,\ldots,K-1~~~(\hbox{where }K\leq M-1).
\end{align*}
%{\color{blue}
%\begin{itemize}
%%\item abscissae and weights
%%\item weighted least squares regression problem
%%\item %\item in the numerical experiments, we use evenly spaced and pdf for weights
%%\item defined the abscissae weights (define d lambda from gautschi)
%%\item can uniquely go from points and weights to three term recurrence coefficients, cite stable lanczos method
%%\item computing inner product in closed form using quadratic over each interval
%\item what you actually have to compute in terms of a matrix equation
%%\item Relation to weighted least squares given in Forsythe paper
%%\item 
%\end{itemize}
%}

%Before proceeding to numerical experiments, w
We briefly comment on the relationship between the spectrum-adapted approximation proposed in this section and the Lanczos approximation 
to $f({\bf A}){\bf b}$, which is given by \cite{druskin}, \cite[Section 13.2]{higham}
\begin{align}\label{Eq:lanc_approx}
{\bf Q}_K f({\bf T}_K){\bf Q}_K^{\top} {\bf b} = ||{\bf b}||_2 {\bf Q}_K f({\bf T}_K){\bf e}_1 ,
\end{align}
where ${\bf Q}_K$ is an $N \times (K+1)$ matrix whose columns form an orthonormal basis for %the Krylov subspace 
%\begin{align*}
${\cal K}_K({\bf A}, {\bf b})=\hbox{span}\left\{{\bf b}, {\bf Ab}, \ldots,{\bf A}^K{\bf b}\right\},$ %the $K$th 
a Krylov subspace.
%\end{align*}
In \eqref{Eq:lanc_approx}, ${\bf T}_K={\bf Q}_K^{\top}{\bf A}{\bf Q}_K$ is a $(K+1) \times (K+1)$ tridiagonal Jacobi matrix. The first column of ${\bf Q}_K$ is equal to $\frac{{\bf b}}{||{\bf b}||}$. The approximation \eqref{Eq:lanc_approx} can also be written as $q_K({\bf A}){\bf b}$, where $q_K$ is the degree $K$ polynomial that interpolates the function $f$ at the $K+1$ eigenvalues of ${\bf T}_{K}$ \cite[Theorem 13.5]{higham}, \cite{saad1992analysis}. Thus, unlike %other 
classical polynomial approximation methods, % such as the truncated Cheybshev expansion, 
the Lanczos method is 
indirectly adapted to the spectrum of ${\bf A}$. 
The Lanczos method differs from proposed method in that  ${\bf T}_{K}$ and the Lanczos approximating polynomial $q_K$ depend on the initial vector ${\bf b}$. %, which is used to initialize the recurrence.
%One difference from the proposed method is that the Lanczos approximating polynomial $q_K$ depends on ${\bf b}$ as well as ${\bf A}$, because the Lanczos process is initialized with ${\bf b}$, and so ${\bf T}_{K}$ depends on ${\bf b}$. 
Specifically, the polynomials $\{\tilde{\pi}_k\}$ generated from the three-term recurrence %of the form \eqref{Eq:three_term} 
\begin{align*}
\gamma_{k+1} \tilde{\pi}_{k+1}(x) = (x-\alpha_k)\tilde{\pi}_k(x)-\gamma_k \tilde{\pi}_{k-1}(x),
\end{align*}
with the $\{\alpha_k\}_{k=0,1,\ldots,K}$ and $\{\gamma_k\}_{k=1,2,\ldots,K}$ coefficients taken from the diagonal and superdiagonal entries of ${\bf T}_K$, respectively, are orthogonal with respect to the piecewise-constant measure
\begin{align*}
\mu(x)=\begin{cases}
0, &x<\lambda_1 \\
\sum_{j=1}^i [\hat{\bf b}(j)]^2, & \lambda_i \leq x < \lambda_{i+1} \\
\sum_{j=1}^N [\hat{\bf b}(j)]^2=1, &\lambda_N \leq x
\end{cases},
\end{align*} 
where $\hat{\bf b}={\bf V}^{\top}{\bf q}_1={\bf V}^{\top}\left(\frac{{\bf b}}{||{\bf b}||}\right)$, and $\hat{\bf b}(j)$ is its $jth$ component \cite[Theorem 4.2]{golub2009matrices}. If %${\bf A}$ has distinct eigenvalues and  
$\hat{\bf b}$ happens to be a constant vector, then $\mu(x)=P_{\lambda}(x)$ from \eqref{Eq:cdf}. %Relate to GFT. 
If ${\bf A}$ is a graph Laplacian, $\hat{\bf b}$ is the graph Fourier transform \cite{shuman2013emerging} of ${\bf b}$, normalized to have unit energy.

\section{NUMERICAL EXAMPLES AND DISCUSSION}

\begin{figure}[bth] 
\begin{minipage}[m]{0.1\linewidth}
~
\end{minipage}
\begin{minipage}[m]{0.44\linewidth}
\centerline{\small{~~~~$|f(\lambda)-p_{10}(\lambda)|$}}
\end{minipage}
\hspace{.01\linewidth}
\begin{minipage}[m]{0.42\linewidth}
\centerline{\small{~~$\frac{||f({\bf A}){\bf b}-p_K({\bf A}){\bf b}||_2^2}{||f({\bf A}){\bf b}||_2^2}$}}
%\centerline{\small{~~$\frac{\sum_{\ell=1}^N \left(f(\lambda_{\ell})-p_K(\lambda_{\ell})\right)^2}{\sum_{\ell=1}^N f(\lambda_{\ell})^2}$}}
\end{minipage}\\
\begin{minipage}[m]{0.1\linewidth}
\centerline{\small{gnp}}
\end{minipage}
\begin{minipage}[m]{0.44\linewidth}
\centerline{~~\includegraphics[width=.95\linewidth]{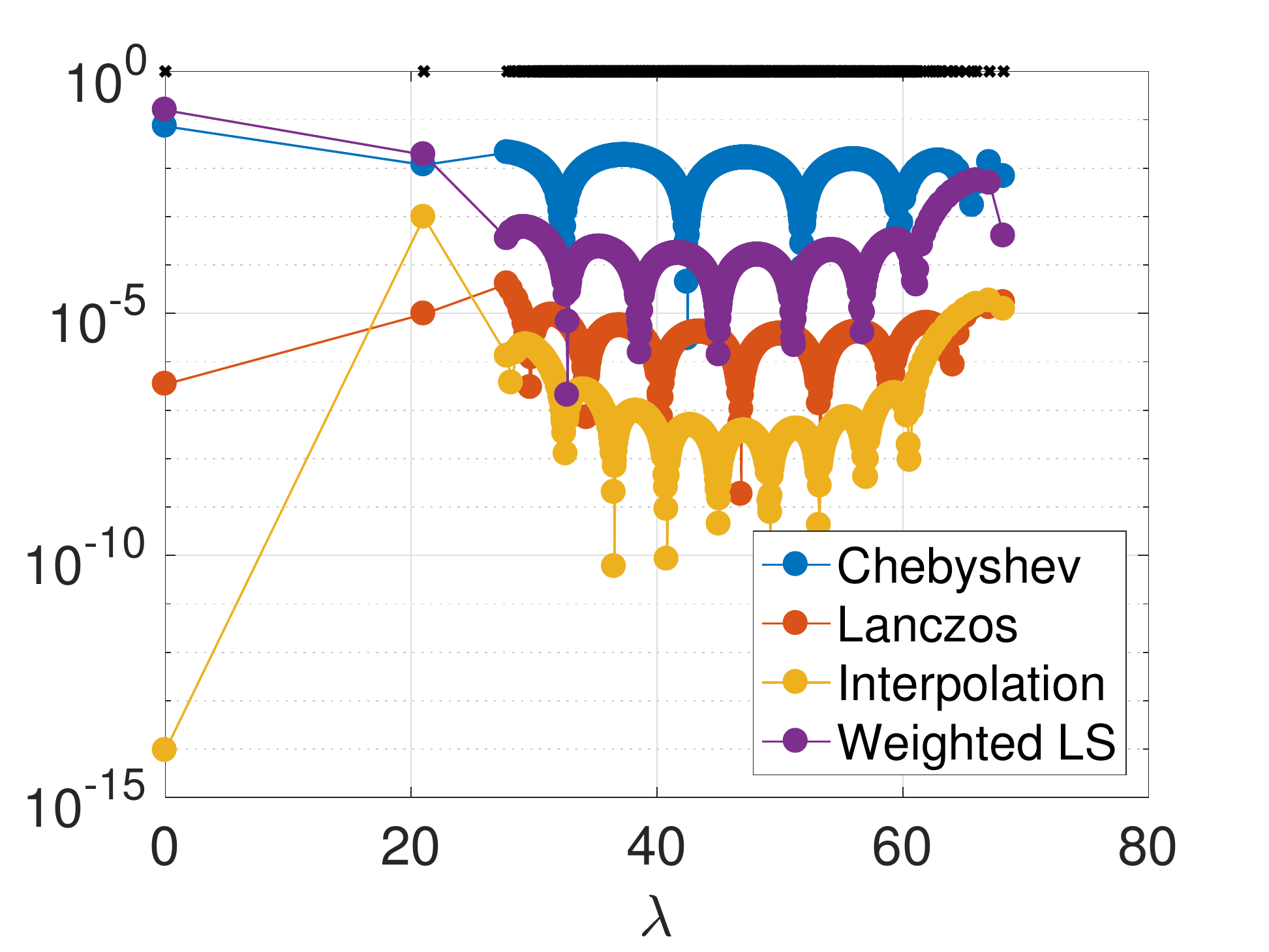}}
\end{minipage}
\begin{minipage}[m]{0.44\linewidth}
\centerline{~~\includegraphics[width=.95\linewidth]{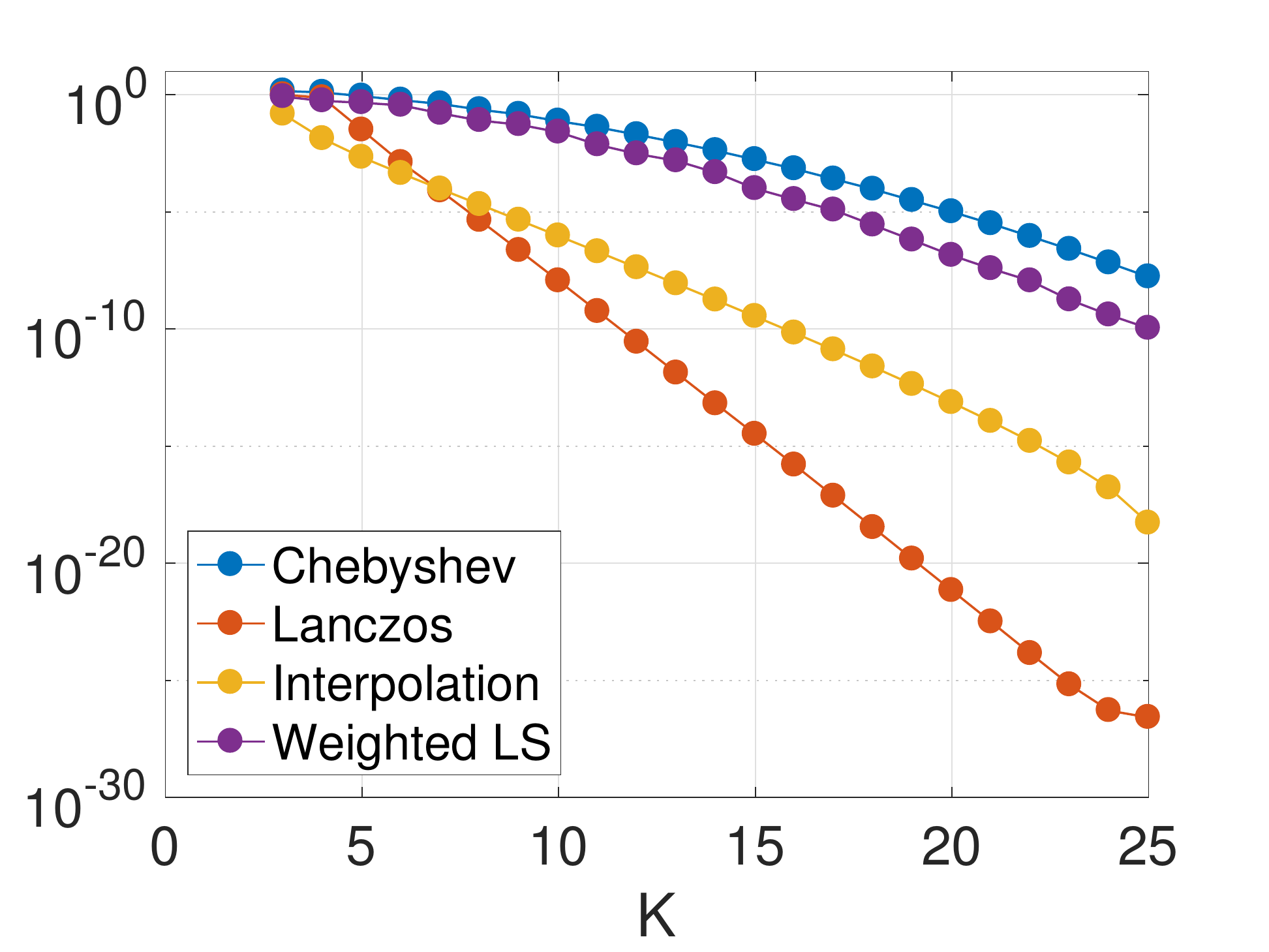}}
\end{minipage}\\
\begin{minipage}[m]{0.1\linewidth}
\centerline{\small{minnesota}}
\end{minipage}
\begin{minipage}[m]{0.44\linewidth}
\centerline{~~\includegraphics[width=.95\linewidth]{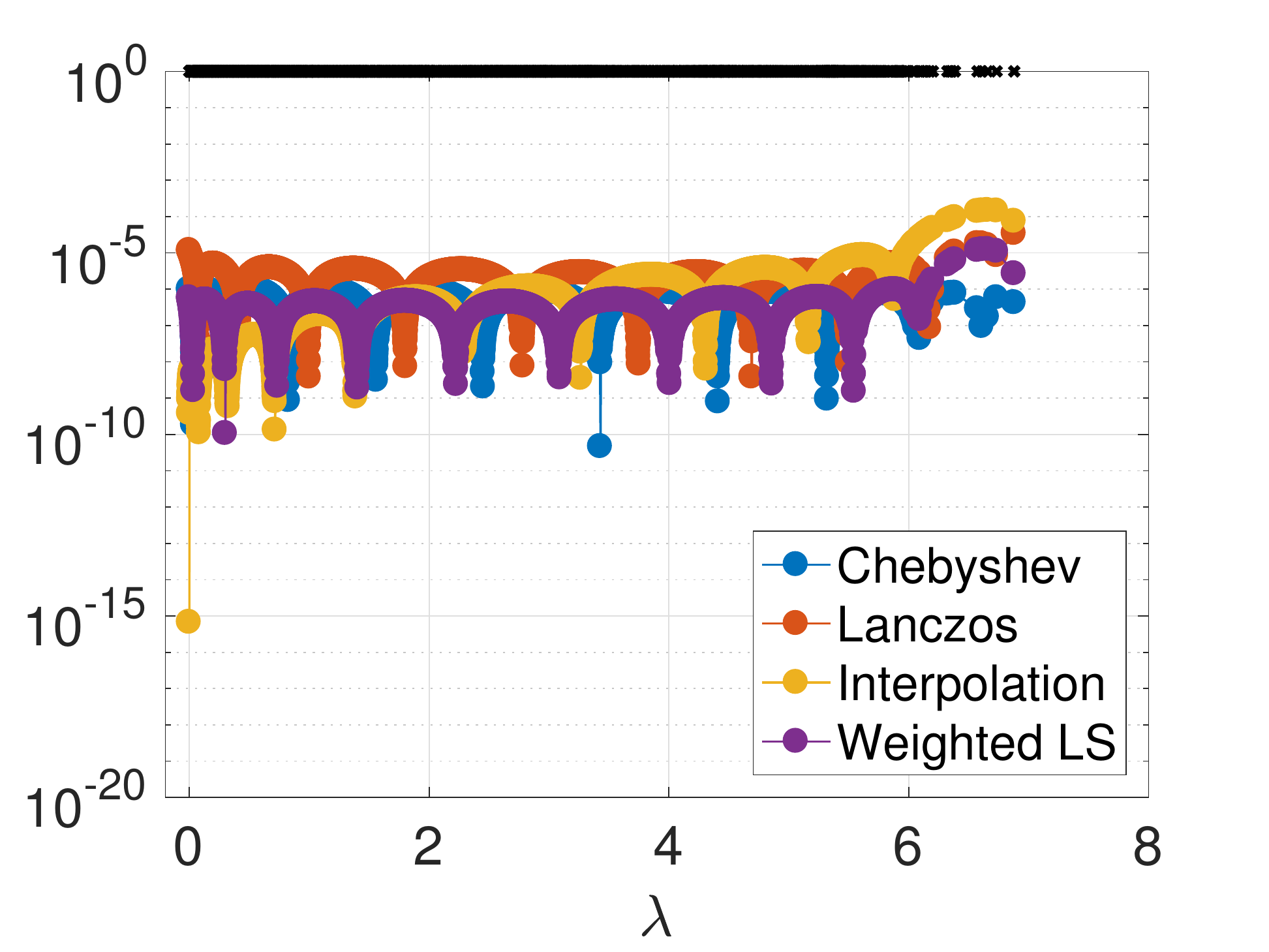}}
\end{minipage}
\begin{minipage}[m]{0.44\linewidth}
\centerline{~~\includegraphics[width=.95\linewidth]{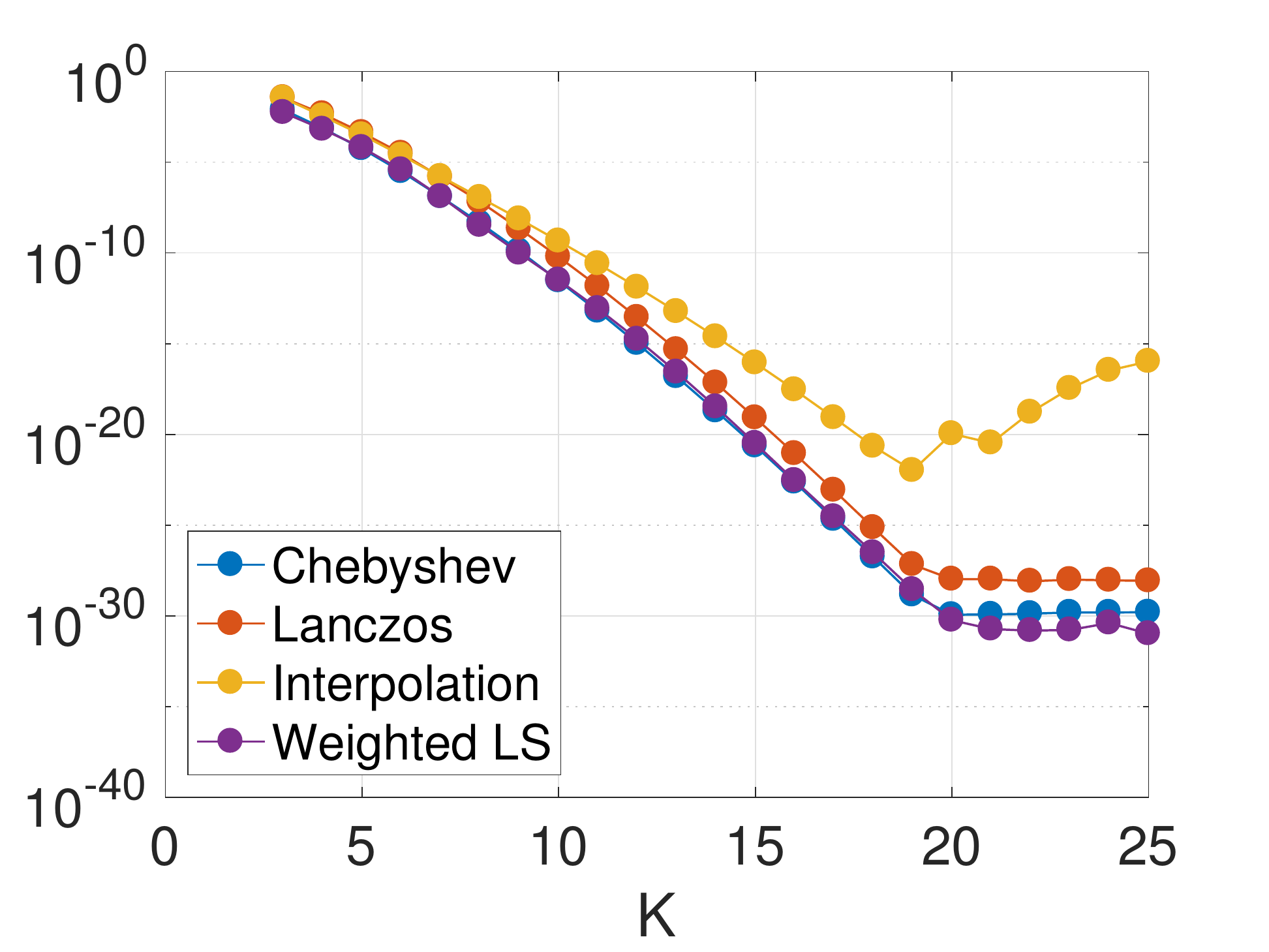}}
\end{minipage}\\
\begin{minipage}[m]{0.1\linewidth}
\centerline{\small{net25}}
\end{minipage}
\begin{minipage}[m]{0.44\linewidth}
\centerline{~~\includegraphics[width=.95\linewidth]{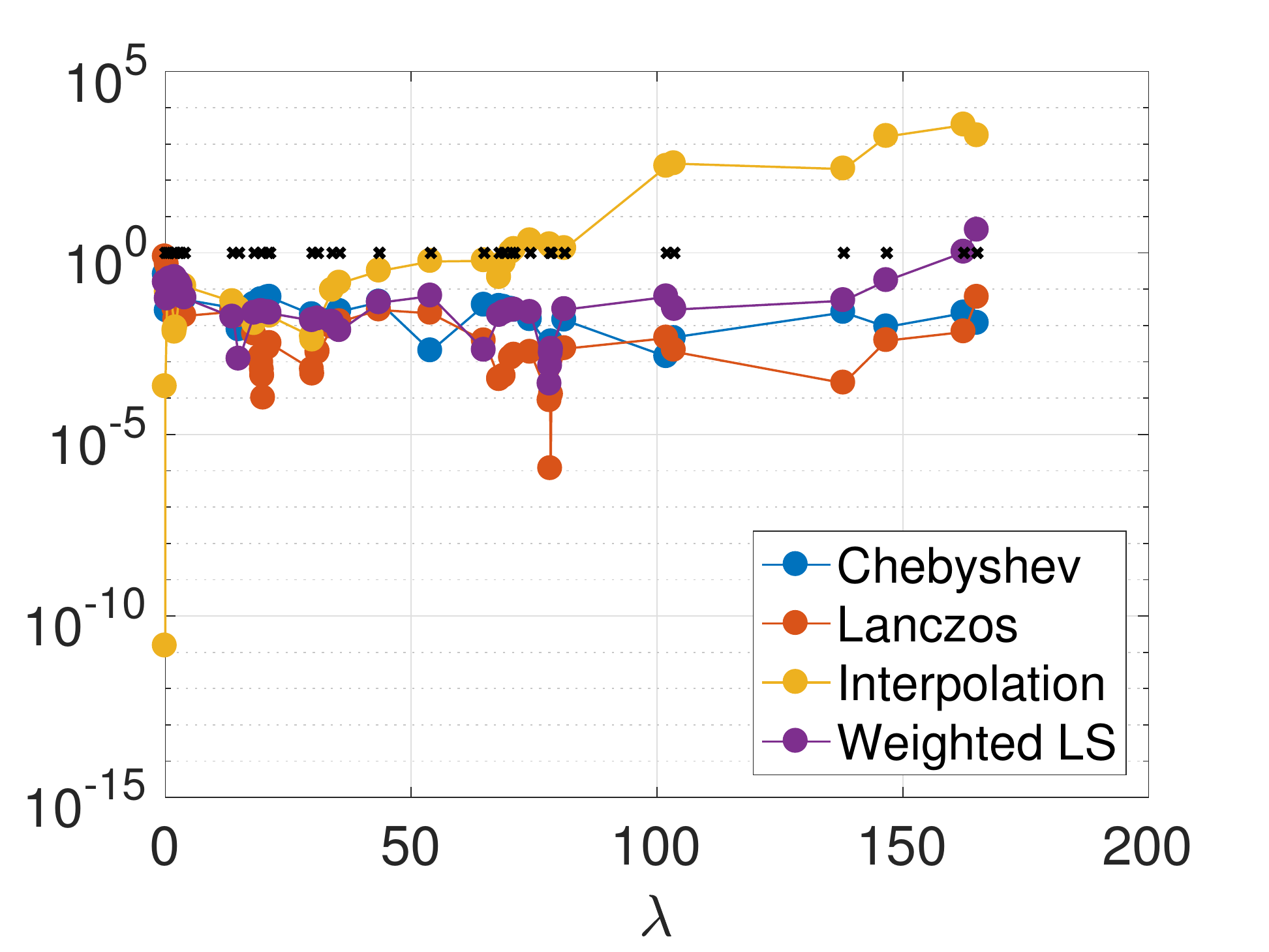}}
\end{minipage}
\begin{minipage}[m]{0.44\linewidth}
\centerline{~~\includegraphics[width=.95\linewidth]{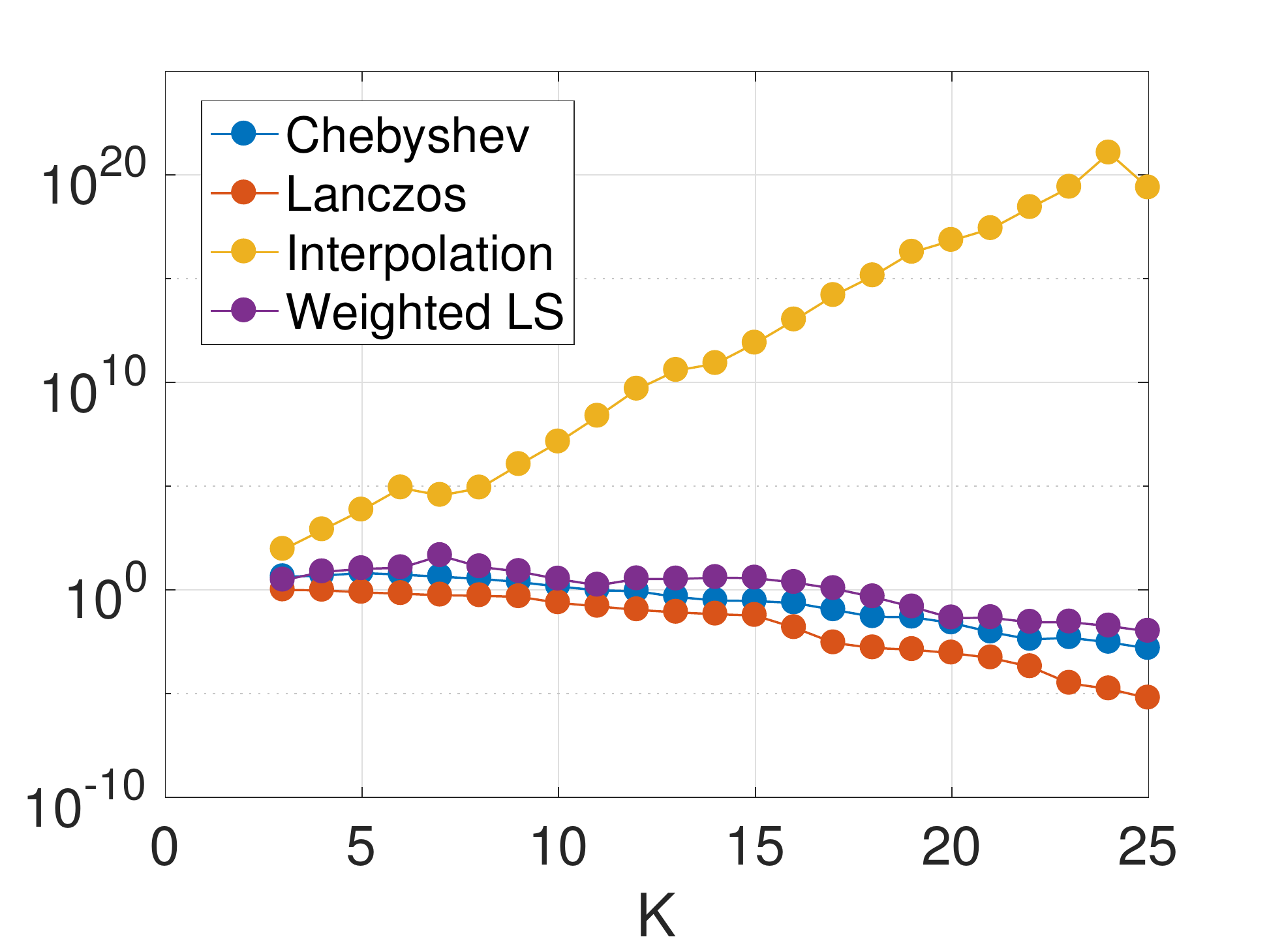}}
\end{minipage}\\
\begin{minipage}[m]{0.1\linewidth}
\centerline{\small{si2}}
\end{minipage}
\begin{minipage}[m]{0.44\linewidth}
\centerline{~~\includegraphics[width=.95\linewidth]{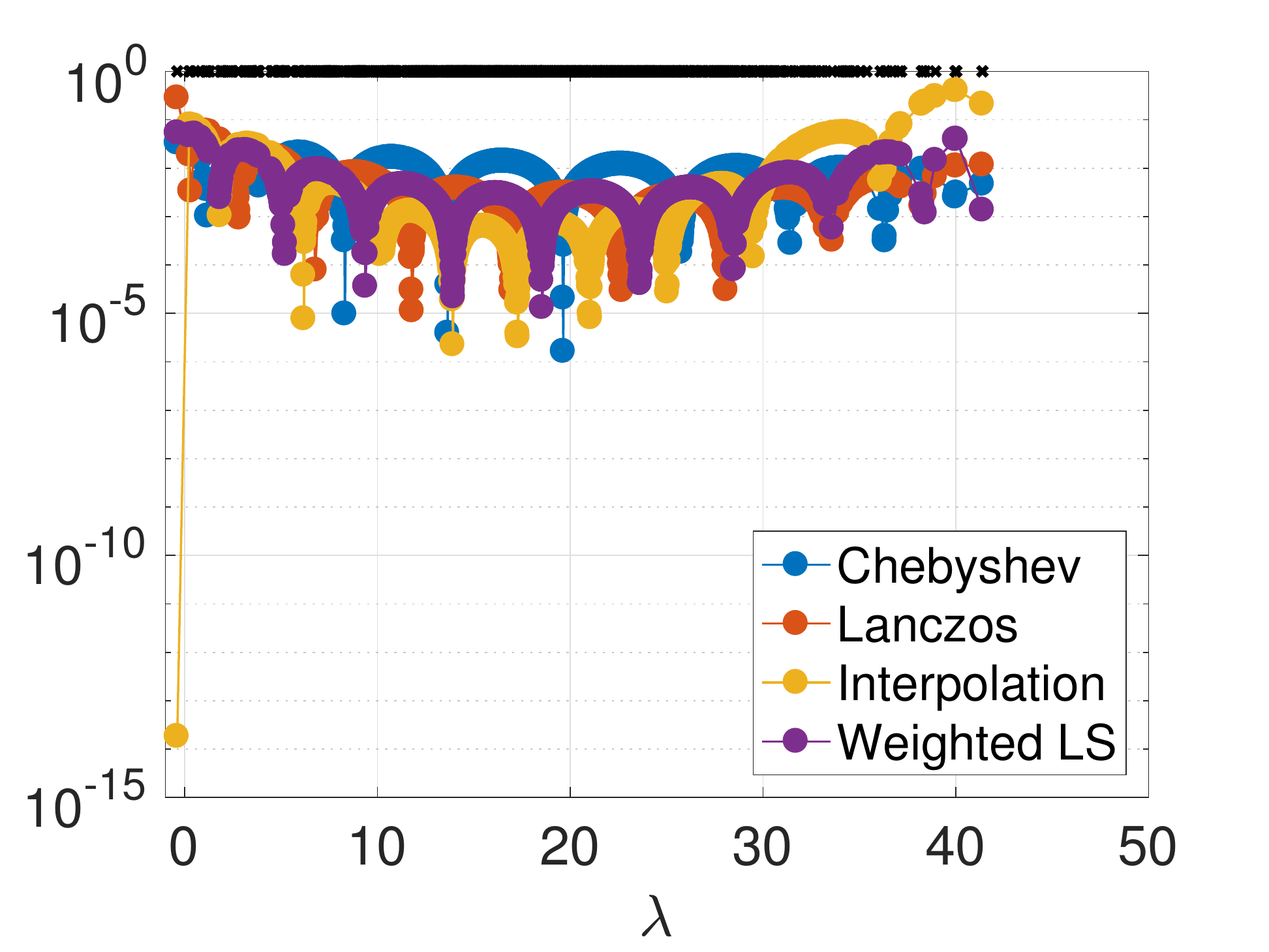}}
\end{minipage}
\begin{minipage}[m]{0.44\linewidth}
\centerline{~~\includegraphics[width=.95\linewidth]{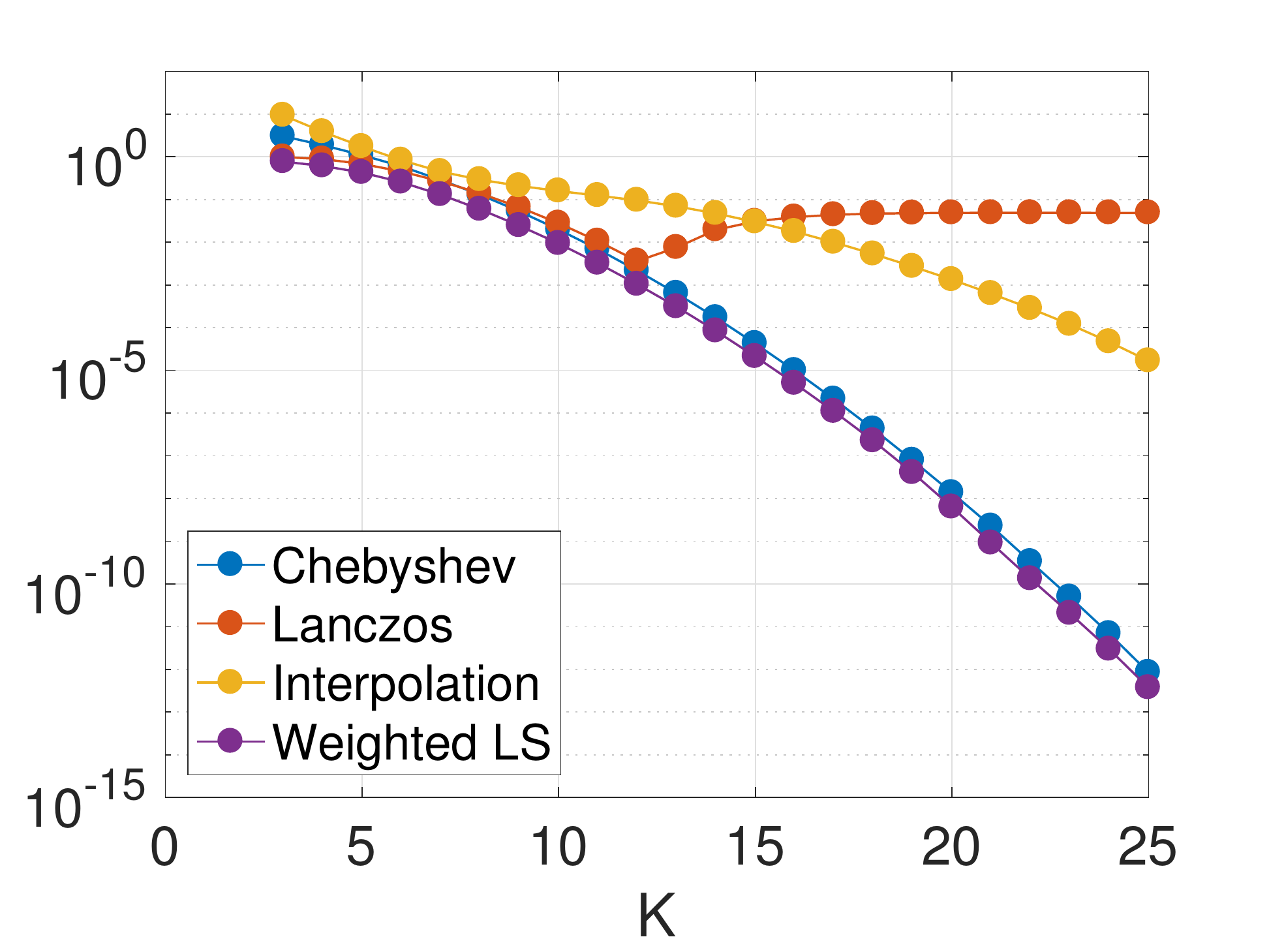}}
\end{minipage}\\
\begin{minipage}[m]{0.1\linewidth}
\centerline{\small{cage9}}
\end{minipage}
\begin{minipage}[m]{0.44\linewidth}
\centerline{~~\includegraphics[width=.95\linewidth]{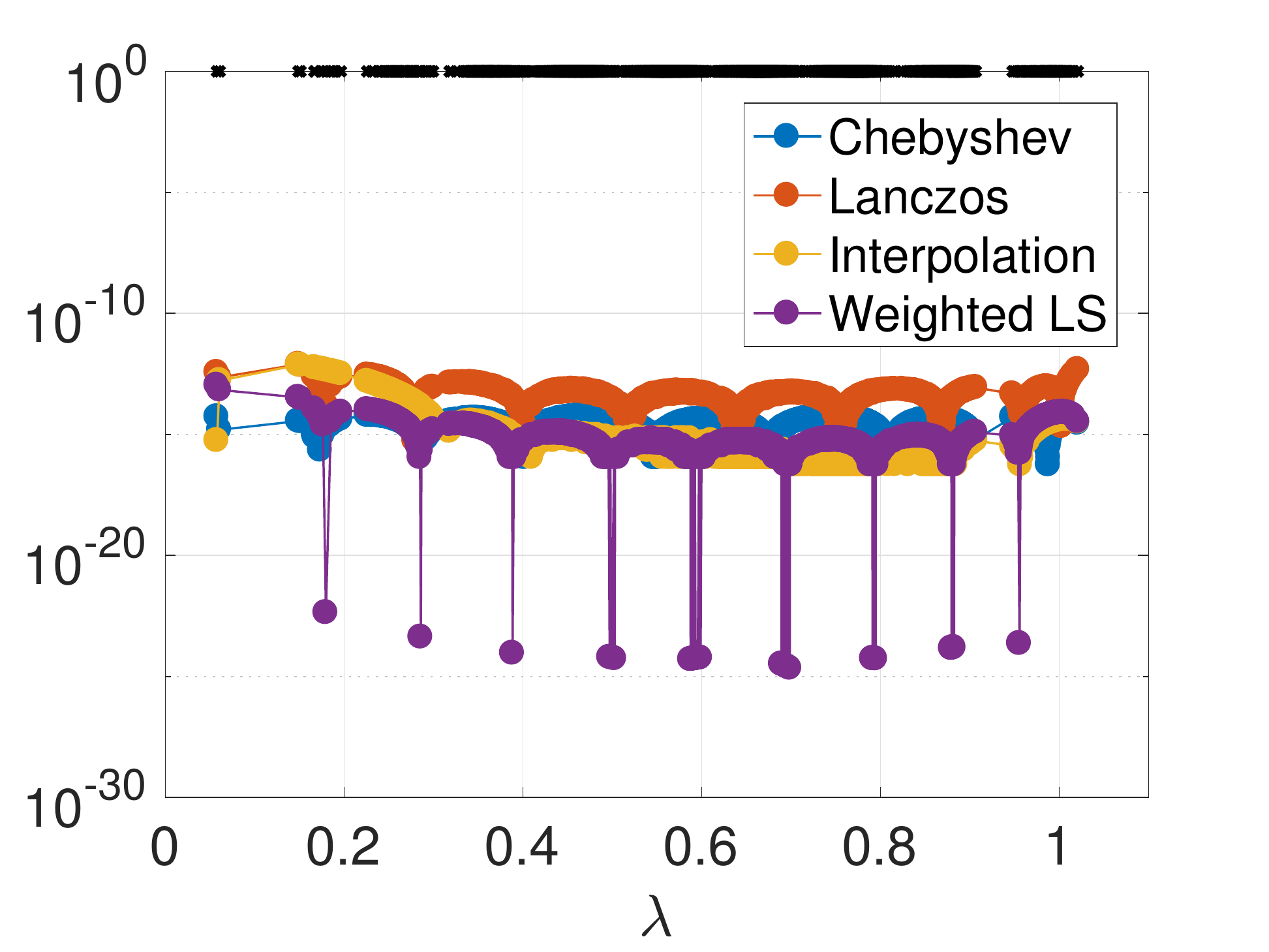}}
\end{minipage}
\begin{minipage}[m]{0.44\linewidth}
\centerline{~~\includegraphics[width=.95\linewidth]{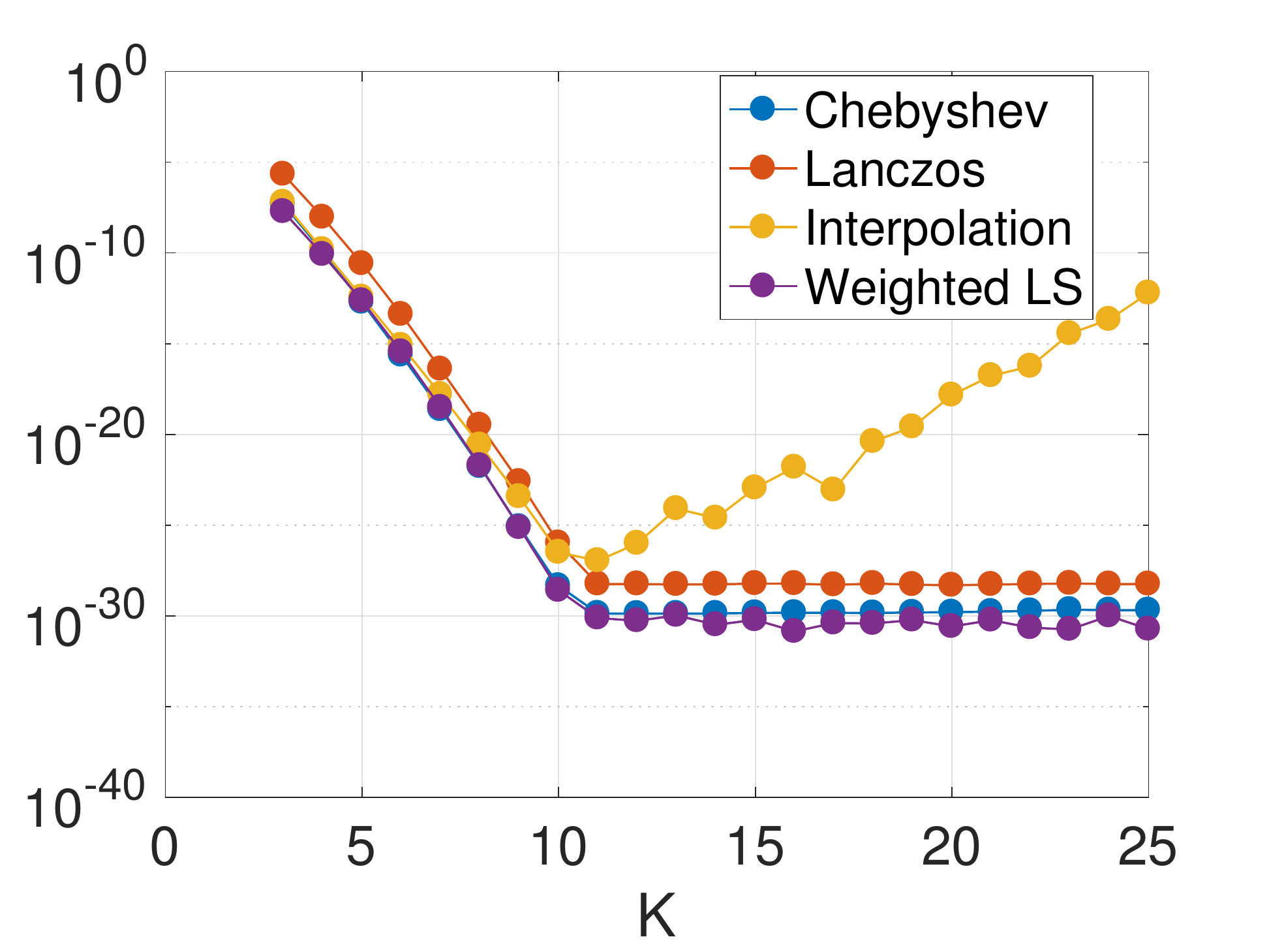}}
\end{minipage}\\
\begin{minipage}[m]{0.1\linewidth}
\centerline{\small{saylr4}}
\end{minipage}
\begin{minipage}[m]{0.44\linewidth}
\centerline{~~\includegraphics[width=.95\linewidth]{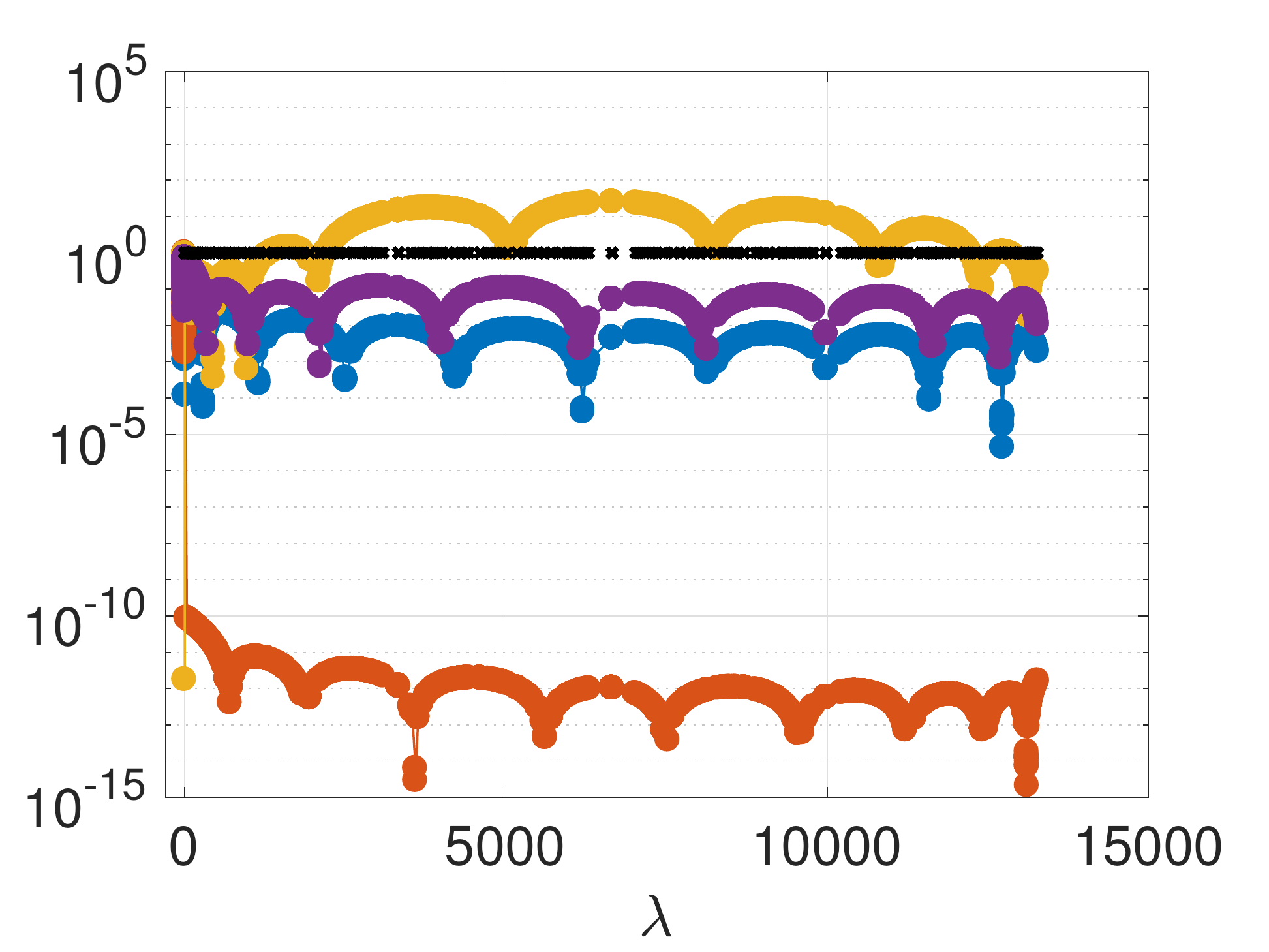}}
\end{minipage}
\begin{minipage}[m]{0.44\linewidth}
\centerline{~~\includegraphics[width=.95\linewidth]{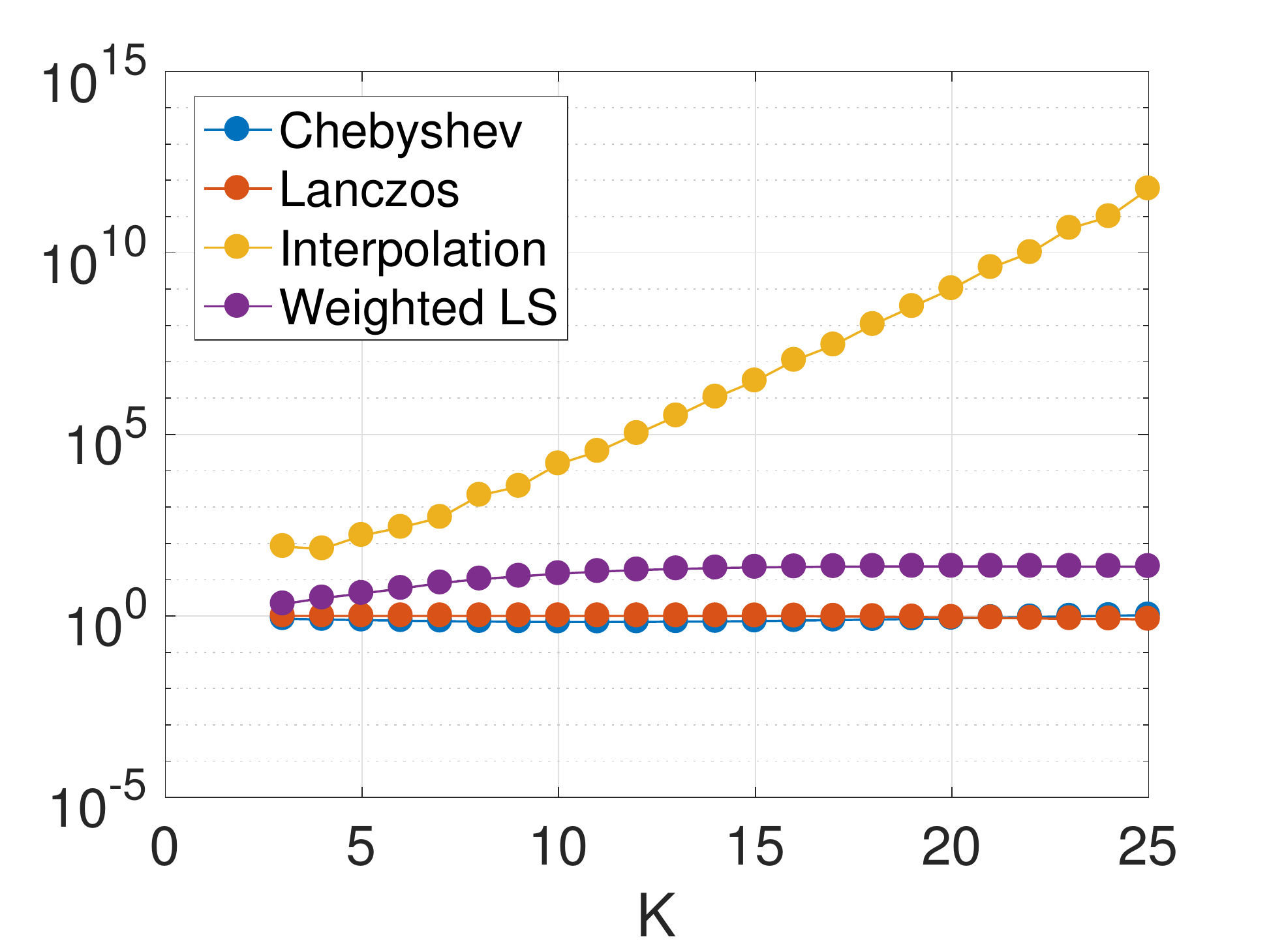}}
\end{minipage}\\
\begin{minipage}[m]{0.1\linewidth}
\centerline{\small{saylr4}}
\centerline{\small{(scaled)}}
\end{minipage}
\begin{minipage}[m]{0.44\linewidth}
\centerline{~~\includegraphics[width=.95\linewidth]{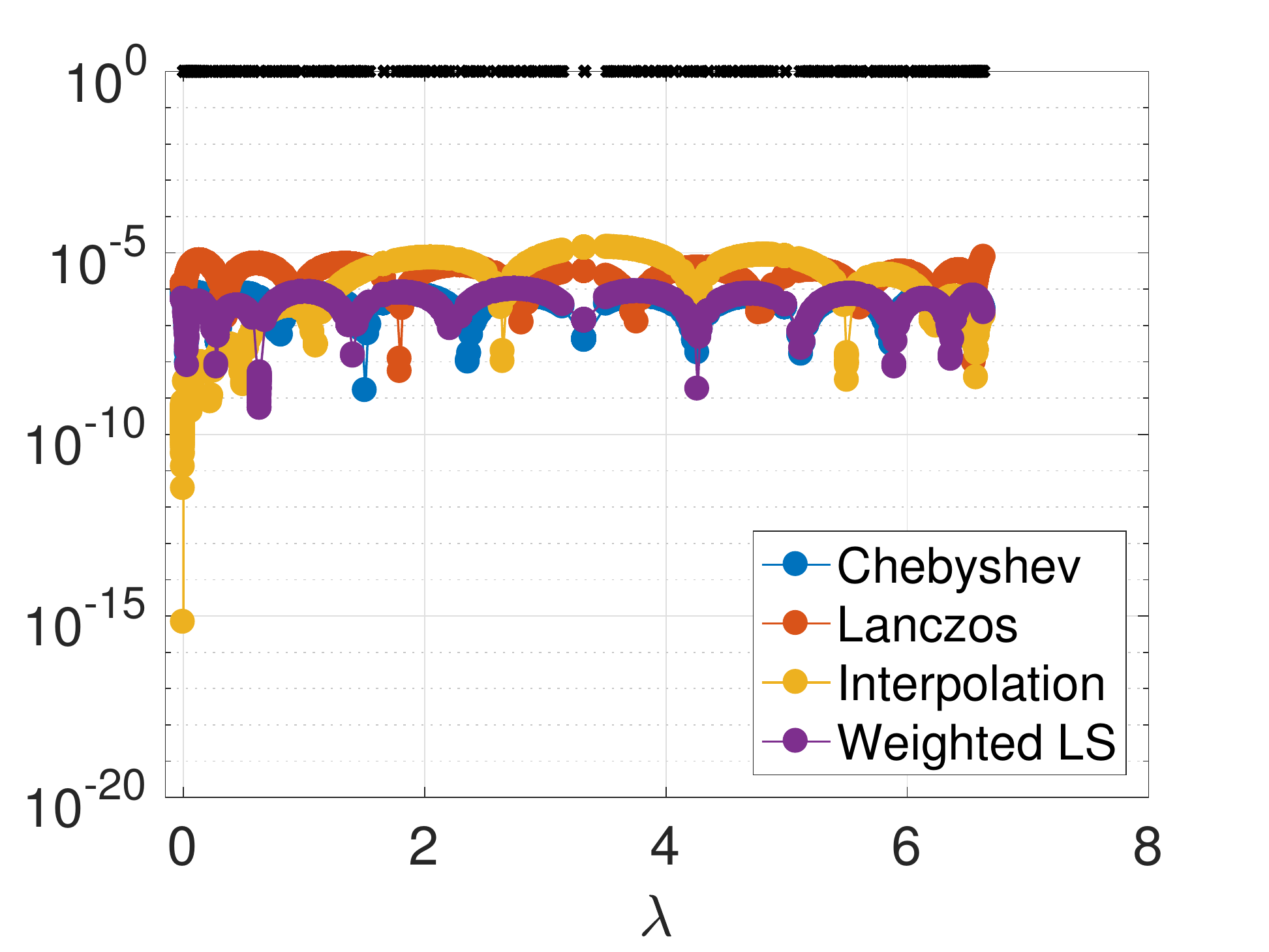}}
\end{minipage}
\begin{minipage}[m]{0.44\linewidth}
\centerline{~~\includegraphics[width=.95\linewidth]{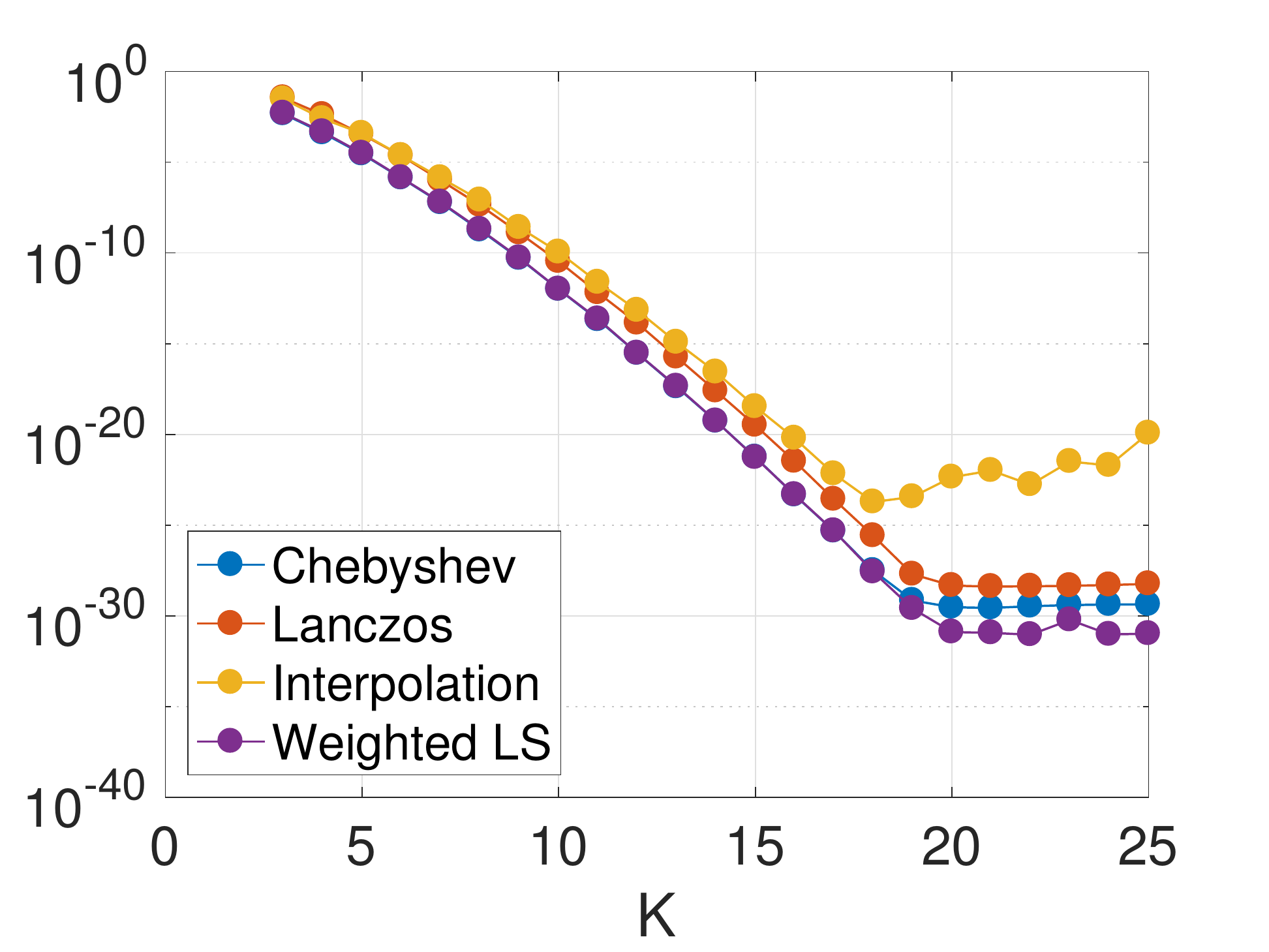}}
\end{minipage}
\caption{Approximations of $f({\bf A}){\bf b}$ with $f(\lambda)=e^{-\lambda}$.}\label{Fig:numerical}
%The first column shows the errors at eigenvalues of {\bf A} given by order 10 polynomial approximations. The second column shows the relative errors in $f({\bf A})$ as the polynomial order increases from 3 to 25.
\end{figure}

We consider the matrix function $f(\lambda)=e^{-\lambda}$, and approximate $f(\bf A){\bf b}$ with ${\bf b}={\bf V1}$ for different matrices ${\bf A}$ and polynomial approximation orders ranging from $K=3$ to 25. %For each matrix ${\bf A}$ and each order, we generate 50 random vectors ${\bf b}$ assuming a uniform distribution in the spectral domain of ${\bf A}$. 
First, we use KPM to estimate the cumulative spectral density function $\tilde{P}_{\lambda}(z)$ with parameters $T=10$,  $J=10$, and $K_{\Theta}=30$, as shown in Fig.\ \ref{Fig:cdf}. Based on the analytical derivative and inverse function of $\tilde{P}_{\lambda}(z)$, we obtain the two proposed spectrum-adapted polynomial approximations for $f(\lambda)$, before computing each $p_K(\bf A){\bf b}$ via the corresponding three-term recursion.
We compare the proposed methods to the truncated Chebyshev expansion and the Lanczos method with the same polynomial order. The results are summarized in Fig. \ref{Fig:numerical}. %The first column shows the errors of each degree 10 approximation at the eigenvalues of ${\bf A}$.
%The approximation errors are measured by the absolute error in \eqref{Eq:exact_error} and the relative error defined as $\sum_{\ell=1}^N \left(f(\lambda_{\ell})-p_K(\lambda_{\ell})\right)^2/\sum_{\ell=1}^N f(\lambda_{\ell})^2$. 
The first column of Fig. \ref{Fig:numerical} displays the %absolute
 errors at all eigenvalues of $\bf A$ for each order 10 polynomial approximation of $f(\lambda)$. The second column examines the convergence of relative errors in approximating $f({\bf A}){\bf b}$ for matrices with various spectral distributions, for each of the four methods. Note that when ${\bf b}$ is a constant vector in the spectral domain of ${\bf A}$, the relative error $\frac{||f({\bf A}){\bf b}-p_K({\bf A}){\bf b}||_2^2}{||f({\bf A}){\bf b}||_2^2}$ is equal to $\frac{\sum_{\ell=1}^N \left(f(\lambda_{\ell})-p_K(\lambda_{\ell})\right)^2}{\sum_{\ell=1}^N f(\lambda_{\ell})^2}$, the numerator of which is the discrete least squares objective mentioned in Section \ref{Se:intro}.
%$\sum_{\ell=1}^N \left(f(\lambda_{\ell})-p_K(\lambda_{\ell})\right)^2$ is equal to $||f({\bf A}){\bf b}-p_K({\bf A}){\bf b}||_2$ when ${\bf b}$ is a constant vector in the spectral domain of {\bf A}, and equivalent to $\mathbb{E}\left[||f({\bf A}){\bf b}-p_K({\bf A}){\bf b}||_2\right]$ if the energy of ${\bf b}$ is uniformly distributed on the spectrum of {\bf A}.

We make a few observations based on the numerical examples: % in Fig. \ref{Fig:numerical}:
\begin{enumerate}
\item The spectrum-adapted interpolation method often works well for low degree approximations ($K\leq 10$), but is not very stable at higher orders due to ill-conditioning.
\item The Lanczos method is more stable than other methods with respect to the width of the spectrum. To demonstrate this, we scaled the saylr4 matrix by multiplying it by the constant $\frac{1}{2000}$ in the bottom row of Fig. \ref{Fig:numerical}. Doing so drastically improves the relative performance of the other methods, even though the spectral distribution is the same.
\item The proposed spectrum-adapted weighted least squares method tends to outperform the Lanczos method for matrices such as si2 and cage9 that have a large number of distinct interior eigenvalues.
%\item Unlike the other methods, the approximating polynomial $p_K(\lambda)$ for the Lanczos method (c.f. first column of Fig. \ref{Fig:numerical}) depends on the choice of ${\bf b}$. 
\item The proposed spectrum-adapted methods, like the Chebyshev approximation, are amenable to efficient distributed computation via communication between neighboring nodes %, as discussed in
 \cite{shuman_distributed_SIPN_2018}. The inner products of the Lanczos method may lead to additional communication expense or severe loss of efficiency in certain distributed computation environments (e.g., GPUs). % is not due to the required inner products. 
\end{enumerate}
%{\color{red}
%%Discussion:
%\begin{itemize}
%%\item Interpolation works alright for low order ($K\leq 10$), but is not very stable at higher orders due to ill-conditioning (community)
%%\item Lanczos is more stable with respect to spectral width (saylr4)
%%\item Some situations in which Lanczos struggles, and our methods may be better
%%\begin{itemize}
%%\item Large number of distinct interior eigenvalues (Si2)
%\item If $b$ does not cover the spectrum, it could be an issue initializing Lanczos with it (add an example)
%\item Speed. Other methods do not require orthogonalization step, should be faster. Check.
%%\item Distributed computations. The other methods do not require the inner products
%%\item Low order approximations?
%%\item Small spectral width?
%%\item Want a single polynomial that does not depend on $b$. Is this ever the case? How does that compare with averaging Lanczos over different $b$?
%%\item Avoid the need for reorthogonalization?
%\end{itemize}
%%\item Look at discrete least squares, which corresponds to constant vector in the spectral domain (${\bf b}={\bf V}{\bf 1}$)
%%\item Different methods to choose abscissae and weights (original saylr4 matrix, different choice of points, same weights, reference paper from Saad) 
%%\end{itemize}
%%Older possibilities:
%%\begin{itemize}
%%\item Lanczos is still the best general purpose method?
%%\end{itemize}
%}

\section{ONGOING WORK} \label{Sec:ongoing}

Our ongoing work includes (i) testing the proposed methods on further applications, such as the estimation of the log-determinant of a large sparse Hermitian matrix; (ii) investigating the theoretical manner and rate at which $p_K({\bf A}){\bf b}$ converges to $f({\bf A}){\bf b}$ for our new approximation methods; (iii) exploring methods to adapt the approximation to the specific matrix function $f$, in addition to the estimated spectral density of ${\bf A}$; (iv) exploring whether these types of approximations can be used to efficiently compute interior eigenvalues in situations where the Lanczos method struggles; and (v) testing whether it is worthwhile to incorporate our proposed methods into the estimation of the eigenvalue counts in \eqref{Eq:eig_count}, in an iterative fashion, since each $\tilde{\Theta}_{\xi_i}({\bf A}){\bf x}^{(j)}$ is itself of the form $f({\bf A}){\bf b}$. 
\clearpage
\section{References}
%\small
\balance
\bibliographystyle{IEEEbib}
{\scriptsize \bibliography{saop_refs}}

\end{document}